%% file: ex_article.tex
\documentclass[hidelinks,onefignum,onetabnum]{siamart220329}
\usepackage{mathdots}
\usepackage{todonotes}
\usepackage{booktabs}
\usepackage{subcaption}
\usepackage{multirow}

\input{ex_shared}

\ifpdf
\hypersetup{
  pdftitle={
  Preconditioned normal equations for solving discretised partial differential equations
  },
  pdfauthor={L. Lazzarino, Y. Nakatsukasa, and U. Zerbinati}
}
\fi

\newtheorem{example}{Example}
\newtheorem{nt}{Notion}
\let\vec\mathbf

\newcommand{\norm}[1]{\left\lVert#1\right\rVert}

\begin{document}

\maketitle

\begin{abstract}
    This paper explores preconditioning the normal equation for non-symmetric square linear systems arising from PDE discretization, focusing on methods like CGNR and LSQR. The concept of {``normal''} preconditioning is introduced, and a strategy to construct preconditioners studying the associated ``normal'' PDE is presented. Numerical experiments on advection diffusion problems demonstrate the effectiveness of this approach in achieving fast and stable convergence.
\end{abstract}

\begin{keywords}
    Preconditioning, normal equation, CGNR, advection-diffusion
\end{keywords}

\begin{MSCcodes}
    65F08, 65N22, 65H10, 65N30, 65N06, 65N55
\end{MSCcodes}
\section{Introduction}
We are interested in studying a linear system of the form
\begin{equation}
    \label{eq:linearsystem}
    A\vec{x} = \vec{b}, \qquad A \in \mathbb{R}^{n \times n}, \vec{x} \in \mathbb{R}^n, \vec{b} \in \mathbb{R}^n,
\end{equation}
where $A$ is invertible,
and often sparse and nonsymmetric, i.e., $A \neq A^T$, and $n\gg 1$. In this case, iterative methods, and in particular Krylov subspace methods, are the most viable option.
Given the non-symmetric nature of the problem, arguably the most common Krylov subspace method is the generalized minimal residual (GMRES) method \cite{Saad}. GMRES is based on minimizing the 2-norm of the residual in increasingly larger Krylov subspaces $\mathcal{K}_k(A,\vec{b}) \coloneq \text{span}(\vec{b},A\vec{b}, \dots , A^k\vec{b})$, where we assume the initial guess $vec{x}_0 = 0$. In the literature, GMRES is often presented as an obvious choice, albeit there is a wide range of applications in which other approaches can have significant advantages.
We here argue for conjugate gradient (CG) like methods \cite{Saad2,Trefethen,Trefethen2,Paige}. Much is known about these methods, such as their convergence behaviour. Hence, we can obtain fast and stable (as we discuss below)
solutions for linear systems of the form \eqref{eq:linearsystem}. To this end, we formulate a linear system equivalent to \eqref{eq:linearsystem}, i.e.~the normal equation
\begin{equation}
    \label{eq:normal}
    B\vec{x} = A^T\vec{b}, \qquad B = A^TA.
\end{equation}
Here, the matrix of coefficients $B$, often called the Gram matrix, is a symmetric positive definite (SPD) matrix. Thus, the system \eqref{eq:normal} can be solved by CG.
The strategy of considering the normal equation and applying CG is well known for solving least-squares problems, and various algorithms have been proposed in this context to implement it. The best known and used are CG normal residual (CGNR) \cite{Saad2}, where the CG method is directly applied to the linear system \eqref{eq:normal}, and LSQR \cite{Paige}, a reformulation of the CGNR algorithm equivalent in exact arithmetic but with better finite precision properties. While alternative methods such as CG normal error (CGNE) \cite{Saad2} exist, that is, applying CG to the dual normal equation
\begin{equation}
    AA^T\vec{y} = \vec{b}, \quad \vec{x} = A^T\vec{y},
\end{equation}
and although we recognize the value in investigating this option, we prioritize CGNR in this work. Our choice is motivated by two key considerations: the availability of \bb{a} mathematically equivalent  implementation, LSQR, which can be made numerically stable and is typically faster, and is widely recognized for its effectiveness in least-squares problems; and the fact that, unlike CGNE, both CGNR and LSQR, at each iteration, minimize the same quantity as GMRES, the standard solver for square systems.  This provides a meaningful basis for comparison with \bb{an} established method.
CGNR and LSQR inherit the known analysis of CG, with the only difference that we now obtain results starting from the matrix $B = A^TA$, instead of $A$. In other words, at each iteration we are minimizing the $A^TA$-norm of the forward error in an increasingly larger Krylov subspace, i.e. at the $k$-th iteration, we have
\begin{equation}
    x_k = \text{arg}\min_{x\in \mathcal{K}_k(A^TA,A^T\vec{b})} \|x - x^*\|_{A^TA},
\end{equation}
where $x^*$ is the exact solution of the linear system. This is equivalent to minimizing the $2$-norm of the residual in $\mathcal{K}_k(A^TA,A^T\vec{b})$. Note that this implies that at each iteration, we are minimizing the same quantity we would minimize in GMRES, but in a different subspace. Moreover, by the convergence analysis of CG, it can be shown that the convergence of CGNR and LSQR is fully determined by the spectral properties of the matrix ${B}$ \cite{Trefethen},
\begin{equation}
    \label{eq:CGbound}
    \frac{\|x^*-x_k\|_B}{\|x^*-x_0\|_B} \leq \underset{p\in P_k}{\bb{\min}} \underset{\lambda \in \Lambda(B)}{\bb{\min}} |p(\lambda)| \leq 2 \left(\frac{\sqrt{\kappa(B)}-1}{\sqrt{\kappa(B)}+1} \right)^k = 2 \left(\frac{\kappa(A)-1}{\kappa(A)+1} \right)^k,
\end{equation}
where $P_k$ is the space of polynomials $p$ of degree less than or equal to $k$ and unit lowest order coefficient, $\Lambda(B)$ is the spectrum of $B$.
The convergence behaviour is not fully characterized by the condition number of $B$, but by its spectral distribution, and that an adaptive bound as in \cite{Strakos2} would be more descriptive of the convergence behaviour of the method. However, bound \eqref{eq:CGbound} indicates that if the condition number of the matrix is very close to 1, then the method will converge very fast independently of the spectral distribution of $B$, and this is enough for our discussion.
The GMRES method is very different from this point of view. Indeed, its convergence is determined by the minimal polynomial of the matrix ${A}$ \cite{Trefethen}, and the spectral properties of $A$ are not enough to fully characterize the convergence of the method. In fact, it can be proven that any nonincreasing GMRES convergence curve is possible for any eigenvalue distribution, \cite{Greenbaum, Strakos}.  Consequently, preconditioning strategies for this method are often derived by heuristic processes. It is worth mentioning that if the matrix $A$ is nearly orthogonally diagonalizable, the convergence of GMRES is well understood in terms of the spectrum of $A$~\cite[Theorem 35.2]{Trefethen}.

\bb{The classical treatment of preconditioning the normal
    equation has been developed primarily for rectangular least-squares
    problems~\cite{Saad2,axelsson1994iterative,scotttuma2025sparse}. However,}
methods like CGNR and LSQR appear to have been sometimes overlooked in the context of non-symmetric square systems mainly because of the need of accessing the action of $A^T$ and because using the Gram matrix squares the condition number of the problem. However, in the context of PDEs discretisation, the action of the adjoint operator can often be computed analytically, or symbolically \cite{dolfinAdjoint}, via integration by parts. Moreover, we argue that by studying the normal equation, we can devise new preconditioners that are competitive in the iteration count in spite of the squared condition number. In addition, as shown in \cite{epperly2025,Maike}, it is possible to obtain a backward stable solution by the preconditioned CGNR or LSQR (when carefully implemented) simply by introducing a small number of iterative refinement steps in the algorithm. In other words, if a good preconditioner is employed, we can find a fast and backward stable solution to the non-symmetric square linear system by the preconditioned CGNR/LSQR method.

In this paper, we first discuss what it means for a preconditioner to be good or ideal in this context, stressing the differences with the GMRES framework.
We revisit a notion of good preconditioner from the least-squares community to introduce the notion of ``normal'' preconditioners for square systems, which characterizes a bigger class of ideal preconditioners.
The choice of preconditioner can be guided by discretisation of a “normal” PDE associated with the continuous problem; more details are provided in Section \ref{sec:funcanal}.
All numerical experiments presented here are performed with CGNR for ease of implementation. We recall that CGNR and LSQR are equivalent in exact arithmetic and, thus, all formal derivations apply to both.
In this paper we mainly focus on CGNR, even though LSQR is sometimes observed to have better numerical stability. This is because we often discuss the Gram matrix $A^TA$ (or its appropriate variant of the form $A^TTA$) directly rather than the matrix $A$, and the two methods are both proven to be stable with iterative refinement~\cite{epperly2025}, provided a careful implementation of the underlying Lanczos iteration~\cite{musco2018stability} is used.
We view this paper as a proof of concept, that solving PDEs  via the preconditioned normal equation, and in general solving non-symmetric linear systems via CG-like methods on the normal equation, can be a viable approach, at least in some cases.


\bb{We use the number of iterations until convergence as a measure of computational effort, for theoretical convenience: it is independent of implementation details and allows us to derive explicit convergence bounds, as in \eqref{eq:CGbound}. We note that more appropriate measures of computational cost should be used for practical comparisons.}

\section{Preconditioning the normal equation}
\label{sec:precond}
In this section, we discuss the use of a preconditioner to improve the convergence of the CGNR method \cite{WathenNormal}.
The most straightforward approach to precondition the normal equation \eqref{eq:normal} is to consider a preconditioned CG method, where the preconditioner is a symmetric positive definite matrix ${G}$.
Such a scheme is equivalent either to solving the systems
\begin{equation}
    \label{eq:normalPrecond}
    {G}^{-1}{B}\,\vec{x} = {G}^{-1}{A}^T\vec{b},\text{ or } {B}\,{G}^{-1}\vec{y} = {A}^T\vec{b}, \text{ with } \vec{x} = {G}^{-1}\vec{y},
\end{equation}
The first system corresponds to a left preconditioned CGNR method, while the second to a right preconditioned CGNR method.
Equation \eqref{eq:CGbound} suggests that, if we can find a preconditioner ${G}$ such that ${G}^{-1}{B}$ or ${B}{G}^{-1}$ have $k$ distinct eigenvalues, then the preconditioned CGNR method will converge in at most $k$ iterations.
In practice, often a preconditioner $G$ is chosen such that ${G}^{-1}{B}$ and ${B}{G}^{-1}$ have clustered eigenvalues.
More often than not, if the matrix ${A}$ arises from the discretisation of a PDE, we have a good physical intuition of what the operator ${A}$ represents.
Using this intuition, a large number of preconditioners have been developed \cite{Elman,Wathen,Benzi,Farrell}.

Assuming that we have a good preconditioner ${P}$ for $A\vec{x} = \vec{b}$, can we use ${G} = {P}^T{P}$ as a preconditioner for the CGNR method?
This question has been addressed by D.~Braess and P.~Peisker \cite{Braess}, and S.~Gratton et al. \cite{Gratton}, and in the recent work by A.~Wathen \cite{WathenNormal}.
First, we would have to agree on what it means for ${P}$ to be a good preconditioner for $A\vec{x} = \vec{b}$. In the case of a symmetric positive definite $A$, and so when CG is available directly on \eqref{eq:linearsystem}, the answer is well known in the literature, also theoretically. However, in the considered case, lacking a general established notion, the following is often used:
\begin{nt}
    \label{nt:classical}
    $P$ is a good preconditioner for  $A\vec{x} = \vec{b}$ if ${P}^{-1}{A}$ has clustered eigenvalues.
\end{nt}
Such preconditioners can be shown to be good for (left-preconditioned) GMRES under some assumptions~\cite{Saad}.
Unfortunately, such notion of goodness does not imply that $G$ is a good preconditioner for CG applied to the normal equation, since a phenomenon called matrix squaring in~\cite{WathenNormal} can cause $G$ to be an arbitrarily bad preconditioner for the normal equation, as shown in \cite{Wathen}.

For this reason, in \cite{WathenNormal} the matrix ${T}:= I - {A}{P}^{-1}$ is considered, and it is observed that
\begin{equation}
    \label{eq:specEquiv}
    {G}^{-1}{B} = {P}^{-1}{P}^{-T}{A}^T{A}\sim {P}^{-T}{A}^T{A}{P}^{-1} = {I} - {T} - {T}^T + {T}^T{T}.
\end{equation}
Since ${G}^{-1}{B}$ is similar to a symmetric positive definite matrix, its eigenvalues are real. Furthermore, since for any symmetric matrix ${S}$ we have $\Lambda({S})\subset [-\norm{{S}}_2, \norm{{S}}_2]$, using Weyl's inequality \cite{horn,stewart}, we can easily see that
\begin{equation}
    1 - 2\norm{{T}}_2 - \norm{{T}}_2^2 \leq \lambda \leq 1 + 2\norm{{T}}_2 + \norm{{T}}_2^2.
\end{equation}
Thus, we see that another possible notion of good preconditioner for CGNR is available. $G$ is a good preconditioner for CGNR if $\norm{I-{A}{P}^{-1}}_2$ is small. \bb{In particular, for $\|T\|_2\leq -1 + \sqrt{2}$ the eigenvalues will be clustered away from zero}. This notion was first introduced by S. Gratton et al. \cite{Gratton}.

We introduce  a notion of goodness for $P$ that implies goodness of $G$ as a preconditioner of CG applied to the normal equation. This will highlight a counterintuitive phenomenon related to the side we decide to precondition \eqref{eq:normal}.
Recall that, from \eqref{eq:CGbound}, the convergence of CGNR is related to the spectrum of $G^{-1}B$ or of $BG^{-1}$. We focus on the left preconditioned normal equation \eqref{eq:normalPrecond}, i.e., on $G^{-1}B$, and observe that from \eqref{eq:specEquiv} we know
\begin{equation}
    \label{eq:crossSquare}
    G^{-1}B \sim (AP^{-1})^T(AP^{-1}).
\end{equation}
It is immediate to see that if $AP^{-1}$ is an orthogonal matrix, then $G^{-1}B$ reduces to the identity and the left preconditioned normal equation in \eqref{eq:normalPrecond} can be solved by means of the CGNR method {in exact arithmetic in exactly one iteration}.
\begin{nt}
    \label{nt:crosseyed}
    $G$ is a good left preconditioner for CG applied to the normal equation in \eqref{eq:normalPrecond} if $AP^{-1}$ has clustered singular values. In this case, we will refer to $P$ as a good normal preconditioner.
\end{nt}
We highlight the ``cross'' effect of the preconditioner action, i.e. to have a good left preconditioner for CGNR, we need a good right preconditioner for $A\vec{x} = \vec{b}$ according to the analogue of Notion 1, but applied to singular values rather than eigenvalues.
Note that this change in the preconditioning side is needed to define this new concept of goodness. Indeed, if $P^{-1}A$ has clustered singular values, $P$ is not necessarily a good preconditioner for CGNR \cite{Wathen}.
The notion of clustering singular values is widely adopted, for example, in the randomised linear algebra community, e.g., Blendenpik~\cite{avron2010blendenpik}, a popular randomised solver for least-squares problems, is based on clustering singular values of $AP^{-1}$. Similarly, the role of clustered singular values in achieving fast convergence of CGNR has been discussed for non-symmetric systems with a Toeplitz coefficient matrix, \cite{chan1993circulant}.

{To make the notion of clustering quantitative, and to address the fact that the eigenvalues of $G^{-1}B$ are the \emph{squares} of the singular values of $AP^{-1}$, recall from \eqref{eq:crossSquare} that $\Lambda(G^{-1}B) = \{\sigma^2 : \sigma\in\Sigma(AP^{-1})\}$. Suppose the singular values of $AP^{-1}$ are clustered in a relative interval $\Sigma(AP^{-1})\subset[1-\epsilon,\,1+\epsilon]$ around one\footnote{Such normalisation can always be achieved by rescaling $P$, since $G=P^TP$ and only the ratio of the extremal singular values enters \eqref{eq:CGbound}}, with $0\le\epsilon<1$. Then $\Lambda(G^{-1}B)\subset[(1-\epsilon)^2,(1+\epsilon)^2]$, so squaring does \emph{not} destroy the clustering, although it weakens it slightly: the eigenvalues remain clustered around $1$, and $\kappa(G^{-1}B)\le\big((1+\epsilon)/(1-\epsilon)\big)^2$. Substituting into \eqref{eq:CGbound}, the per-iteration CGNR contraction factor is bounded by
\[
    \frac{\sqrt{\kappa(G^{-1}B)}-1}{\sqrt{\kappa(G^{-1}B)}+1} \le \frac{(1+\epsilon)-(1-\epsilon)}{(1+\epsilon)+(1-\epsilon)} = \epsilon .
\]
In other words, clustering the singular values of $AP^{-1}$ within a relative radius $\epsilon$ of $1$ guarantees an error reduction of at least a factor $\epsilon$ per iteration. The squaring caused by the normal equation only turns the relative radius $\epsilon$ into a condition number $\big((1+\epsilon)/(1-\epsilon)\big)^2$, whose square root recovers exactly the factor $(1+\epsilon)/(1-\epsilon)$. This is the precise sense in which $AP^{-1}$ having clustered singular values makes $P$ a good normal preconditioner.}

We illustrate the cross effect of the preconditioner with the one-dimensional advection diffusion problem, with a finite difference discretisation. We remark that example \ref{ex:cross-eyed} is meant to be illustrative of the concept of normal preconditioning and the ``cross'' effect of the preconditioner action, rather than offering a practical preconditioning strategy (which will be discussed in Section \ref{sec:funcanal}). To this aim, we construct dense left and right preconditioners using direct methods that result in an optimal number of iterations, when the ``cross'' effect is taken into account.
\begin{example}[{Normal preconditioning}]
    \label{ex:cross-eyed}
    We consider the classical advection diffusion ODE in one dimension, i.e.
    \begin{equation}
        \label{eq:advectionDiffusion1D}
        -\nu  u^{\prime\prime} + \beta u^\prime = f \text{ in }(a,b)\subset\mathbb{R}, \text{ and } u(a)=u_a,\; u(b)=u_b,\; \nu,\beta \in \mathbb{R}_{\geq 0}.
    \end{equation}
    We consider neither diffusion nor advection dominated regimes, i.e. $\nu\approx \beta$, and a centered finite difference discretisation over an equi-spaced mesh of step-size $h$.
    Such discretisation results in the linear system
    \begin{equation}
        \label{eq:linsysAdvDiff1D}
        A\vec{x} = \vec{b}, \qquad A=\text{tridiag}\left(-\frac{\nu}{h^2}- \frac{\beta}{2h}, \frac{2\nu}{h^2}, -\frac{\nu}{h^2}+\frac{\beta}{2h}\right),\; \vec{b}_j=f(a+jh).
    \end{equation}
    {Here the unknowns $\vec{x} = (x_1,\dots,x_n)^T$ are the nodal values at the interior points $a+jh$, $j=1,\dots,n$. The Dirichlet data $u(a)=u_a$ and $u(b)=u_b$ enter the right-hand side through the first and last equations, where the boundary nodes $x_0 = u_a$ and $x_{n+1} = u_b$ are known and are moved to the right-hand side, so that
    \begin{equation}
        \label{eq:linsysAdvDiff1Drhs}
        \vec{b}_1 = f(a+h) + \Big(\tfrac{\nu}{h^2}+\tfrac{\beta}{2h}\Big) u_a, \qquad
        \vec{b}_n = f(a+nh) + \Big(\tfrac{\nu}{h^2}-\tfrac{\beta}{2h}\Big) u_b,
    \end{equation}
    while $\vec{b}_j = f(a+jh)$ for $2\le j\le n-1$. In this way the boundary values are encoded in $\vec{b}$, as they must be for the discrete solution to satisfy the prescribed boundary conditions.}
    We are now interested in solving \eqref{eq:linsysAdvDiff1D} via the left preconditioned CGNR method.
    Notice that, since the matrix $A$ is tridiagonal, we could solve the linear system \eqref{eq:linsysAdvDiff1D} in $O(n)$ operations, yet we are here interested in the performance of CGNR, for illustrative purposes.
    \bb{We compare} different choices of preconditioner $P$ and always choose $G=P^TP$.
    In particular, since we are interested in highlighting the cross phenomenon, we consider as $P$ the $R$-factor in both the QR and RQ factorisation of $A$, and the left and right polar factors $(A^TA)^{1/2}$ and $(AA^T)^{1/2}$.
    The results for different meshes of $(0,1)$ with step-size $h$, expressed in terms of degrees of freedom $n$, are shown in Table \ref{tab:cross-eyed} for $\nu=\beta=1$ and $u_a=0$, $u_b=1$.
\begin{table}
        \centering
        \caption{Comparison of the number of iterations for different preconditioners for the left preconditioned normal equation in \eqref{eq:normalPrecond}.
            The CGNR method was terminated when the \bb{relative residual norm $\lVert Ax - b \rVert / \lVert b \rVert$ of the original system $A\vec{x}=\vec{b}$, i.e.~the unpreconditioned residual,} was less than $10^{-10}$.
            If the method did not converge in $1000$ iterations, we marked the number of iterations with a dash.
        }
        {
            \label{tab:cross-eyed}
            \begin{tabular}{c|c|c|c|c}
                \toprule
                $n$  & QR & RQ & $Q(A^TA)^{1/2}$ & $(AA^T)^{1/2}Q$ \\
                \midrule
                10   & 1  & 12 & 1               & 3               \\
                100  & 1  & -  & 1               & 3               \\
                1000 & 1  & -  & \bb{1}          & \bb{3}          \\
                \bottomrule
            \end{tabular}
        }
    \end{table}
\end{example}
The previous example highlights another crucial difference between preconditioning \eqref{eq:normal} and \eqref{eq:linearsystem}, i.e.~the ideal preconditioner for the normal equation is not unique.
In fact, while clustering the eigenvalues of $AP^{-1}$ to a single point, in the ideal scenario, results in $P$ being $A^{-1}$ up to scaling, the same is not true for the normal equation.
Clustering singular values of $AP^{-1}$ to a single point, in the ideal scenario, amounts to $AP^{-1}$ being an orthogonal matrix up to scaling.
The space of orthogonal matrices is a manifold of dimension $n(n-1)/2$~\cite{edel98, epperly2025}. Hence, the space of ideal normal preconditioners for the normal equation is much larger than the space of ideal preconditioners for the linear system.
To demonstrate this, in the previous example, we constructed two different ideal preconditioners for the normal equation, i.e.~the $R$-factor in the QR factorisation of $A$ and the polar factor $(A^TA)^{1/2}$.

\section{Functional Analytic Perspective}
\label{sec:funcanal}
We explore the ideas presented above in this section, with a focus on singularly perturbed PDEs {\cite{Verhulst}}, the prototypical example being the advection diffusion equation in the advection dominated regime.
Consider the two-dimensional advection diffusion equation, i.e.~
\begin{equation}
    \label{eq:advectionDiffusion2D}
    \mathcal{L}u := -\nu \Delta u + \vec{\beta}\cdot \nabla u = f \text{ in }\Omega\subset\mathbb{R}^d, \text{ and } u = \bb{0} \text{ on }\partial\Omega, \text{ with } \nu\ll \norm{\beta},\;\nabla\cdot \beta =0.
\end{equation}
In particular, we will focus on a finite element discretisation.
To this end, we fix a discrete space $V_h\subset H^1_0(\Omega)$ and consider the weak formulation of \eqref{eq:advectionDiffusion2D} in the discrete space $V_h$, i.e.~find $u_h\in V_h$ such that
\begin{equation}
    \label{eq:weakForm}
    (\hat{\mathcal{L}}u_h,v_h)\coloneqq\nu (\nabla u_h, \nabla v_h)_{L^2(\Omega)} + (\beta\cdot \nabla u_h,v_h)_{L^2(\Omega)} = (f,v_h)_{L^2(\Omega)}\text{ for any } v_h\in V_h.
\end{equation}
Given a fixed basis $\{\varphi_i\}_{i=1}^{n}$ for $V_h$, the weak formulation introduced above can be expressed as a linear system:
\begin{equation}
    A\vec{x}=\vec{b}, \text{ with } A_{ij}=(\hat{\mathcal{L}}\varphi_j,\varphi_i)\text{ and } b_j = (f,\varphi_j)_{L^2(\Omega)}.
\end{equation}
From the previous expression, we can work out what PDE is represented by $A^T$, in fact
\begin{equation}
    A^T_{ij} = A_{ji} = (\hat{\mathcal{L}}\varphi_i,\varphi_j) = (\varphi_i, \hat{\mathcal{L}}^*\varphi_j) = (\hat{\mathcal{L}}^*\varphi_j,\varphi_i),
\end{equation}

where $\hat{\mathcal{L}}^*:H^1(\Omega)\subset L^2(\Omega) \to H^1(\Omega)\subset L^2(\Omega)$ is the restriction of the Hilbert adjoint of $\hat{\mathcal{L}}:H^1(\Omega)\to L^2(\Omega)$ on $H^1(\Omega)$ extended to the all $L^2(\Omega)$.
Thus, this is the Hilbert adjoint computed with respect to the $L^2(\Omega)$, i.e.
\begin{equation}
    (\hat{\mathcal{L}}u,v)_{L^2(\Omega)} =  (u,\hat{\mathcal{L}}^*v)_{L^2(\Omega)} \text{ for any } u,v\in H^1(\Omega).
\end{equation}
Notice that this last equation is simply saying that $\hat{\mathcal{L}}^*$ is the adjoint of $\hat{\mathcal{L}}$ computed by means of integration by parts.
\bb{Boundary conditions play an essential role in this adjoint calculation. The integration by parts that defines $\hat{\mathcal{L}}^*$ generates the boundary contributions
\begin{equation}
\nu\int_{\partial\Omega}(v\,\partial_n u - u\,\partial_n v)\,ds - \int_{\partial\Omega}(\beta\cdot n)\,u v\,ds,
\end{equation}
which all vanish because $u,v\in V_h\subset H^1_0(\Omega)$ satisfy the homogeneous Dirichlet condition imposed in \eqref{eq:advectionDiffusion2D}. Consequently $\hat{\mathcal{L}}^* = -\nu\Delta - \beta\cdot\nabla$, again equipped with homogeneous Dirichlet boundary conditions is well defined, since $\nabla\cdot\beta = 0$.
}

    {Let us denote with $V_h^\prime$ the dual space of $V_h$}. It is important to notice that $A:V_h\to V_h^\prime$ and $A^T:V_h\to V_h^\prime$, in fact, the matrix $A^T$ is neither the Hilbert adjoint of $A$ nor the Banach adjoint of $A$, but it is the Galerkin approximation of the Hilbert adjoint of the operator discretised by $A$.

Due to the range of $A$ being $V_h^\prime$ while the domain of $A^T$  is $V_h$, we observe that even if the normal equation \eqref{eq:normal} makes sense from the finite dimensional linear algebra point of view, it does not have any meaning from the functional analysis point of view, in that it is performing operations on inconsistent function spaces.
We remark that this inconsistency is not an issue if one considers the primal variational formulation of the problem \eqref{eq:weakForm}, but is rather introduced by our choice of solving the linear system via the normal equation.
For this reason, we need to be careful when considering the normal equation in the context of the discretisation of partial differential equations.
Let us consider a Riesz map $\tau:V_h^\prime\to V_h$ discretised by the matrix $T$.
To remove the inconsistency mentioned above, we thus consider the following modified normal equation,
\begin{equation}
    \label{eq:normalPDE}
    A^T TA\vec{x} = A^T T\vec{b},
\end{equation}
where now $A^TTA:V_h\to V_h$ and $A^T T: V_h^\prime \to V_h$.
\begin{remark}
    Given a fixed scalar product $(\cdot,\cdot)$ compatible with the space $V_h$, i.e.~such that there exists a Hilbert space $\big(V,(\cdot,\cdot))$ such that $V_h\subset V$ and $(\cdot,\cdot)|_{V_h\times V_h} = (\cdot,\cdot)$, the Riesz map $\tau$ is elliptic, and therefore the operator $T$ is symmetric positive definite, \cite{Kirby}.
    Since the operator $T$ is symmetric positive definite, we can consider the Cholesky factorisation of $T$, i.e.~$T = C^TC$, and rewrite \eqref{eq:normalPDE} as
    \begin{equation}
        \label{eq:normalPDECh}
        (CA)^T(CA)\vec{x} = (CA)^TC\vec{b}.
    \end{equation}
    Thus, \eqref{eq:normalPDE} is the normal equation associated with the linear system
    \begin{equation}
        \label{eq:linSysPDE}
        CA\vec{x}= C\vec{b}.
    \end{equation}
\end{remark}

Analogously to \eqref{eq:specEquiv} , we can show the following spectral equivalence
\begin{equation}
    \label{eq:specEquiv2}
    \begin{aligned}
        (P^TTP)^{-1}A^TTA & = P^{-1}T^{-1}P^{-T}A^TTA       \\
                          & \sim T^{-1}P^{-T}A^TATP^{-1}    \\
                          & =  T^{-1}(AP^{-1})^TT(AP^{-1}).
    \end{aligned}
\end{equation}
Therefore, an ideal preconditioner is such that $(AP^{-1})^TT(AP^{-1}) = T$, i.e. $AP^{-1}$ is orthogonal with respect to the inner product induced by $T$. More generally, as stated by the following theorem, this implies that the convergence behaviour of preconditioned CG applied to \eqref{eq:normalPDE} is governed by the singular values of $AP^{-1}$ with respect to the inner product induced by $T$, or equivalently, the singular values of the matrix $CAP^{-1}$.

\bb{Weighting the normal equation by an SPD matrix, i.e. considering
    $A^T T A$ in place of $A^T A$, is a familiar device in the
    least-squares community, where the weight typically encodes
    measurement covariance or is chosen purely to precondition~%
    \cite{scotttuma2025sparse,bjorck2024numerical}; here, by contrast,
    $T$ is not a free parameter but is forced by the requirement that~\eqref{eq:normalPDE}
    is well posed as a map between the correct function spaces,
    which is what determines the associated normal PDE.}

For any matrix $M$, it is possible to generalize the standard SVD \cite{horn} to the inner product induced by $T$. We define the SVD so that the $T$-singular values of $M$, denoted by $\sigma_{T,i}(M)$, are the square roots of the eigenvalues of the Gram matrix $T^{-1}M^TTM$. This leads to the decomposition $M = U \Sigma V^T T$, with $U^T T U = I$ and $V^T T V = I$, and \bb{$\Sigma$ a non-negative real} diagonal matrix. We refer to the set of $T$-singular values of $M$ as $\Sigma_T(M)$. The largest and smallest $T$-singular values are denoted by $\sigma_{T,\max}(M)$ and $\sigma_{T,\min}(M)$, respectively, and the $T$-condition number is defined as their ratio:
$\kappa_T(M) = \sigma_{T,\max}(M)/\sigma_{T,\min}(M)$.
\begin{remark}
    \b
    The SVD introduced above can be equivalently written as $M=U \Sigma V^T$, with $U^TTU = I$ and $V^TT^{-1}V = I$. That is, we can rewrite the definition of the $T$-SVD with a more standard notation, avoiding including $T$ in the definition, but imposing a notion of "$T^{-1}$- orthogonality" on $V$.
\end{remark}
\begin{theorem}
    \label{thm:conv}
    The relative error at the $k$-th iteration of the CG method applied to the preconditioned equation
    \begin{equation}
        \label{eq:PrecnormalPDE}
        G_T^{-1}B_T\vec{x} = G_T^{-1}A^TT\vec{b}, \quad G_T = P^TTP, \quad B_T = A^TTA,
    \end{equation}
    with exact solution $x^*$, is such that
    \begin{equation}
        \frac{\|x^* - x_k \|_{B_T}}{\|x^* - x_0 \|_{B_T}} \leq \underset{p\in P_k}{\bb{\min} \underset{\,\,\sigma \in \Sigma_T(AP^{-1})}{\bb{\max}}} |p(\sigma^2)| \leq 2 \left( \frac{\kappa_T(AP^{-1}) -1}{\kappa_T(AP^{-1}) +1} \right)^k,
    \end{equation}
    where $P_k$ is the space of polynomials $p$ of degree less than or equal to $k$ and unit lowest order coefficient.
\end{theorem}

\begin{proof}
    Recall that since $T$ is symmetric positive definite, it always exists the matrix $T^{1/2}$, that is as well symmetric positive definite. Observe that, \bb{left multiplying by}  $T^{1/2}P$, we can rewrite equation~\eqref{eq:PrecnormalPDE} as:
    \begin{equation}
        \label{eq:PrecnormalPDE-rewritten}
        T^{-1/2}P^{-T}A^T T A P^{-1} T^{-1/2} y = T^{-1/2}P^{-T}A^T Tb, \quad x=P^{-1}T^{-1/2}y.
    \end{equation}
    To simplify the notation, let us define the matrix $Z := P^{-T}A^T T A P^{-1}$.
    Further, denoting by $y^*$ the exact solution of equation~\eqref{eq:PrecnormalPDE-rewritten}, the following relation holds
    \begin{equation}
        \begin{aligned}
            \|y^*- y\|_{T^{-1/2}Z T^{-1/2}} & = (y^*-y)^TT^{-1/2}Z T^{-1/2}(y^*-y)                         \\
                                            & = (x^*-x)^T(P^T T^{1/2}) T^{-1/2}Z T^{-1/2}(T^{1/2}P)(x-x^*) \\
                                            & = (x^*-x)^T P^T Z P(x-x^*)                                   \\
                                            & = (x^*-x)^T A^TTA(x-x^*)                                     \\
                                            & = \|x^*-x\|_{A^TTA}.
        \end{aligned}
    \end{equation}
    Moreover, by the definition of $T$-singular values given above, and by the similarity transformation
    \begin{equation}
        T^{-1}Z = T^{-1/2}(T^{-1/2}ZT^{-1/2})T^{1/2},
    \end{equation}
    we have that
    \begin{equation}
        \label{eq:singVsEige}
        \sigma^2_{T,i}\left(AP^{-1}\right) = \lambda_i\left(T^{-1}Z\right) = \lambda_i \left( T^{-1/2}ZT^{-1/2} \right).
    \end{equation}
    Thus,
    \begin{equation}
        \begin{aligned}
            \frac{\|x^* - x_k \|_{B_T}}{\|x^* - x_0 \|_{B_T}} & = \frac{\|y^* - y_k\|_{T^{-1/2}ZT^{-1/2}}}{\|y^* - y_0\|_{T^{-1/2}ZT^{-1/2}}}                                     \\
                                                              & \leq \underset{p\in P_k}{\bb{\min}} \underset{\,\,\lambda \in \Lambda(T^{-1/2}ZT^{-1/2})}{\bb{\min}} |p(\lambda)| \\
                                                              & = \underset{p\in P_k}{\bb{\min}} \underset{\,\,\sigma \in \Sigma_T(AP^{-1})}{\bb{\min}} |p(\sigma^2)|
        \end{aligned}
    \end{equation}

    Finally, combining \eqref{eq:CGbound} and \eqref{eq:singVsEige}, we obtain
    \begin{equation}
        \begin{aligned}
            \frac{\|x^* - x_k \|_{B_T}}{\|x^* - x_0 \|_{B_T}} & \leq 2\left( \frac{\sqrt{\kappa(T^{-1/2}ZT^{-1/2})} -1}{\sqrt{\kappa(T^{-1/2}ZT^{-1/2})} +1} \right)^k \\
                                                              & = 2\left( \frac{\sqrt{\kappa(T^{-1}Z)} -1}{\sqrt{\kappa(T^{-1}Z)} +1} \right)^k
        \end{aligned}
    \end{equation}
\end{proof}
The role of different inner products in determining the convergence of preconditioned minimal residual methods has been investigated in \cite{PestanaWathen}.
\begin{remark}
    When applying the CG method to equation \eqref{eq:normalPDE} , we are minimizing the residual in the $T$-norm. Indeed,
    \begin{equation}
        \|x^* -x_k\|_{A^TTA} = (x^*-x_k)^TA^T T A(x^* - x_k) = \|A(x^* - x_k)\|_T = \|r_k\|_T.
    \end{equation}
    Thus, at each iteration, CG on \eqref{eq:normalPDE} is minimizing the same quantity as GMRES on \eqref{eq:linearsystem}, but in a different Krylov subspace and a different norm.
\end{remark}
\begin{remark}
    We have already pointed out that the choice of ideal preconditioner for the normal equation is not unique, and is much larger than the choice of ideal preconditioner for the original linear system, due to the goal being orthogonality of $AP^{-1}$, with respect to the inner product induced  $T$, rather than reducing $AP^{-1}$ to an identity.
    The same observation also holds for the normal equation \eqref{eq:normalPDE}. Furthermore, once we have chosen a Riesz map $\tau$ , any absolutely continuous measure with respect to the Lebesgue measure will induce another valid Riesz map.
    Therefore, varying the choice of the Riesz map $\tau$ , the possible choices of ideal preconditioners for the normal equation become even richer.
\end{remark}
\begin{remark}
    It is worth noticing that $T$ is a dense matrix. Hence, for an efficient implementation of the CGNR method, it is imperative to use a matrix-free implementation of $A^TTA$.
\end{remark}
\begin{remark}
    \label{eq:remarkRiesz}
    The appropriate choice of the Riesz map is problem dependent. It is worth noting that, in our discussion, the Riesz map was introduced to address the inconsistency between function spaces that arises in the formation of the normal equation. Consequently, within our framework, the Riesz map does not entirely define the preconditioner. The relationship between the selection of a Riesz map, an inner product, and a preconditioner has been examined in \cite{gunnel2014note}.
    In our context, the choice of the Riesz map $\tau$ determines the ``normal'' PDE associated with the normal equation \eqref{eq:normalPDE} and not every choice of $\tau$ is consistent with the function spaces involved in the formulation of the continuous problem, thus leading to meaningless ``normal'' PDEs.
    For example, the choice of the Riesz map induced by the $L^2$ inner product, even if formally can be associated with an anisotropic diffusion equation as ``normal'' PDE, will not be consistent with the function spaces of the continuous problem.
    Hence, we suggest choosing a consistent Riesz map $\tau$ and then computing the normal PDE as discussed in Example \ref{ex:H1InnerProduct}.
\end{remark}

Based on what has been described so far, we can outline a general framework that looks at the normal equation and that can potentially be used to construct practical preconditioners for various PDEs, in particular for singularly perturbed PDEs.
This framework proceeds by first identifying the appropriate class of Riesz maps $\tau:V' \to V$ corresponding to the function spaces underlying the continuous PDE.
This choice will determine the ``normal'' PDE associated with the normal equation, whose discretisation $B$ will guide the choice of a suitable normal preconditioner $G$ for the normal equation.
In Example \ref{ex:H1InnerProduct}, we illustrate how, using this framework, we can precondition the advection diffusion equation in the advection dominated regime via the discretisation of the reaction diffusion operator.

\begin{example}[Preconditioning via reaction diffusion]
    \label{ex:H1InnerProduct}
    Once again, we consider the weak formulation of the advection diffusion equation in the advection dominated regime, i.e.~\eqref{eq:weakForm}.
    In particular, we have here adopted a $\mathcal{P}^1$ finite element discretisation on a uniform triangulation, obtained by splitting each square of a uniform grid into two triangles cutting along the south-west to north-east diagonal.
    As Riesz map $\tau$, we consider the one associated with the $H^1_0(\Omega)$ inner product, i.e.
    \begin{equation}
        \label{eq:RieszH1}
        (\nabla \tau(f_h),\nabla v_h)_{L^2(\Omega)} = \nu^{-1}\langle f_h, v_h\rangle \text{ for any } v_h \in V_h.
    \end{equation}
    Once again, we consider the continuous analogue of the weak formulation \eqref{eq:weakForm} and \eqref{eq:RieszH1}, and observe that $TA$ is the Galerkin approximation of the following weak formulation
    \begin{align}
        \nu(\nabla w,\nabla v)_{L^2(\Omega)} & = \nu (\nabla u, \nabla v)_{L^2(\Omega)} + (\beta\cdot \nabla u,v)_{L^2(\Omega)} \\
                                             & = (\nu \nabla u - \beta u,\nabla v)_{L^2(\Omega)},
    \end{align}
    where the last equality follows from integration by parts and the fact that $u$ and $v$ vanish on the boundary of $\Omega$.
    Hence, we can conclude that $\nabla w = \nabla u - \nu^{-1}\Pi_{\nabla} (\beta u)$, where $\Pi_{\nabla}$ is the $L^2(\Omega)$ projection onto the space of gradients of functions in $H^1_0(\Omega)$.
    Using the previous result, we can compute the weak formulation approximated by $A^TTA$, i.e.
    \begin{align}
        \label{eq:weakFormNormalH1}
        (\hat{\mathcal{L}}^*w,v)_{L^2(\Omega)} & = \nu(\nabla w,\nabla v)_{L^2(\Omega)} - (\nabla\cdot(\beta w),v)_{L^2(\Omega)}                                  \\
                                               & = \nu ( \nabla u, \nabla v)_{L^2(\Omega)} + \nu^{-1}(\Pi_{\nabla}(\beta u),\beta v)_{L^2(\Omega)}                \\
                                               & - \nu^{-1} (\Pi_{\nabla}(\beta u),\nabla v)_{L^2(\Omega)} - \nu^{-1}(\nabla u ,\beta v)_{L^2(\Omega)}            \\
                                               & = \nu ( \nabla u, \nabla v)_{L^2(\Omega)} + \nu^{-1}(\Pi_{\nabla}(\beta u),\Pi_{\nabla}(\beta v))_{L^2(\Omega)}.
    \end{align}
    Notice that we have used the fact that $(\Pi_{\nabla}(\beta u),\nabla v)_{L^2(\Omega)} = -(\nabla u,\Pi_{\nabla}(\beta v))_{L^2(\Omega)}$ and that $(\Pi_{\nabla}(\beta u),\beta v)_{L^2(\Omega)} = (\Pi_{\nabla}(\beta u),\Pi_{\nabla}\beta v)_{L^2(\Omega)}$.
    Therefore, we can precondition the normal equation \eqref{eq:normalPDE} using a discretisation of the variational problem: find $u_h\in V_h$ such that
    \begin{equation}
        \label{eq:discreteReactionDiffusion}
        \nu ( \nabla u_h, \nabla v_h)_{L^2(\Omega)} + \nu^{-1}(\Pi_{\nabla}(\beta u_h),\Pi_{\nabla}(\beta v_h))_{L^2(\Omega)} = \langle g_h,v_h \rangle, \text{ for any } v_h\in V_h,
    \end{equation}
    where $g_h\in V_h^\prime$ represents the right-hand side of \eqref{eq:normalPDE}.

        \begin{table}
        \centering
	\caption{\bb{Number of iterations of the LU-preconditioned CGNR method (exact inversion of \eqref{eq:discreteReactionDiffusion}) for different values of the viscosity $\nu$, the streamline-diffusion parameter $\delta$, and the mesh size. The wind is fixed to $\beta = (1,0)$, $f \equiv 1$, and the method is terminated when the unpreconditioned residual norm is below $10^{-5}$. With $\delta = 0$ the stabilisation is switched off, $A^T T A$ coincides with \eqref{eq:discreteReactionDiffusion} and CGNR converges in a single iteration; for a fixed $\delta$ the count grows with $\delta$ but is independent of the mesh. The last row of each block uses the standard SUPG value $\delta = h/2|\beta|$.}}
        \label{tab:reactionDiffusionLU}
        \begin{tabular}{c|c|c|c|c}
            \toprule
            $\nu$ & $\delta$ & $32\times 32$ & $64\times 64$ & $128\times 128$ \\
            \midrule
            $1\cdot 10^{-2}$    & $0$          & 1 & 1 & 1 \\
                                & $10^{-4}$    & 2 & 2 & 2 \\
                                & $10^{-3}$    & 4 & 4 & 4 \\
                                & $h/2|\beta|$ & 17 & 10 & 7 \\
            \midrule
            $5\cdot 10^{-3}$    & $0$          & 1 & 1 & 1 \\
                                & $10^{-4}$    & 3 & 3 & 3 \\
                                & $10^{-3}$    & 5 & 5 & 5 \\
                                & $h/2|\beta|$ & 30 & 17 & 10 \\
            \midrule
            $2.5\cdot 10^{-3}$  & $0$          & 1 & 1 & 1 \\
                                & $10^{-4}$    & 3 & 3 & 3 \\
                                & $10^{-3}$    & 3 & 3 & 3 \\
                                & $h/2|\beta|$ & 53 & 30 & 17 \\
            \midrule
            $1.25\cdot 10^{-3}$ & $0$          & 1 & 1 & 1 \\
                                & $10^{-4}$    & 3 & 3 & 3 \\
                                & $10^{-3}$    & 8 & 8 & 8 \\
                                & $h/2|\beta|$ & 100 & 61 & 31 \\
            \bottomrule
        \end{tabular}
    \end{table}
    In Table \ref{tab:reactionDiffusionLU} we present the iteration count for a preconditioned CGNR method, where as preconditioner we used an exact inversion of \eqref{eq:discreteReactionDiffusion} via an LU factorisation.
    \bb{To keep the discretisation well resolved on the coarser meshes we stabilise it with the standard streamline-upwind (SUPG) value $\delta_K = h_K/2|\beta|$. As Figure~\ref{fig:layers} shows, this choice resolves the boundary and internal layers cleanly, without the spurious oscillations that a smaller stabilisation leaves on under-resolved meshes.}
    Since we adopted a streamline diffusion, we do not expect CGNR to converge in a single iteration.
    \bb{To make this precise, Table~\ref{tab:reactionDiffusionLU} reports the iteration count over a range of $\delta$. When $\delta = 0$ the stabilisation is switched off, the discrete operator $A^T T A$ coincides with the discretisation \eqref{eq:discreteReactionDiffusion} of the normal PDE, and its exact LU inversion yields convergence in a single iteration for every value of $\nu$ and every mesh size. As $\delta$ is increased the streamline-diffusion term perturbs $A$ away from the operator that \eqref{eq:discreteReactionDiffusion} discretises, so the mismatch between $A^T T A$ and the preconditioner grows: at the small fixed values $\delta = 10^{-4}, 10^{-3}$ this costs only a few iterations and remains independent of the mesh. The standard SUPG value $\delta = h/2|\beta|$ is markedly larger, already exceeding $\nu$ on the meshes considered, so the mismatch is severe, and because $\delta$ now scales with $h$ the count both grows sharply as $\nu\to0$ and loses its mesh independence, halving under each refinement.}
    \bb{This is not a defect of the stabilisation but of using the \emph{unstabilised} projected operator \eqref{eq:discreteReactionDiffusion} to precondition a stabilised $A$. The remedy, adopted for the practical sparse preconditioner below, is simply to stabilise the preconditioner consistently, i.e.~to include the same streamline term in the surrogate. With that correction the iteration counts remain robust in $\nu$ and mesh independent even at the SUPG value, as reported in Tables~\ref{tab:reactionDiffusionAMG}--\ref{tab:reactionDiffusionAMG3}.}
    Unfortunately, the mass matrix $\nu^{-1}(\Pi_{\nabla}(\beta u_h),\Pi_{\nabla}(\beta v_h))_{L^2(\Omega)}$ is singular, hence as $\nu$ vanishes it becomes harder and harder to invert it via a multigrid method.
    For this reason, we opted to use as preconditioner the discretisation of the reaction diffusion operator corresponding to the variational problem: find $u_h\in V_h$ such that
    \begin{equation}
        \label{eq:discreteReactionDiffusion2}
        \begin{aligned}
            \nu ( \nabla u_h, \nabla v_h)_{L^2(\Omega)} &{}+ \nu^{-1}(\beta u_h,\beta v_h)_{L^2(\Omega)} \\
            &\bb{{}+ \delta(\beta\cdot\nabla u_h,\beta\cdot\nabla v_h) = \langle g_h,v_h \rangle}, \text{ for any } v_h\in V_h,
        \end{aligned}
    \end{equation}
    \bb{with the same right-hand side $g_h\in V_h^\prime$ as in \eqref{eq:discreteReactionDiffusion},}
    and inverting the preconditioner using PETSc GAMG \cite{Adams}, a smoothed aggregation algebraic multigrid implementation.
    \bb{It is worth clarifying that, even for a constant wind $\beta$ and in the singularly perturbed limit $\nu\to 0$, the dominant reaction term of \eqref{eq:discreteReactionDiffusion2} is a scaled \emph{mass} matrix and not a multiple of the identity. Indeed, for constant $\beta$ one has $\nu^{-1}(\beta u_h,\beta v_h)_{L^2(\Omega)} = \nu^{-1}\lvert\beta\rvert^2 (u_h, v_h)_{L^2(\Omega)}$, i.e.~a rescaling of the $\mathcal{P}^1$ mass matrix $M$. On any fixed mesh $M$ is not spectrally clustered since its eigenvalues span a fixed $\mathcal{O}(1)$ range independent of $\nu$ and so inverting \eqref{eq:discreteReactionDiffusion2} by algebraic multigrid is a genuine computational task and the preconditioner is far from trivial. This is confirmed by the non-trivial iteration counts of Table~\ref{tab:reactionDiffusionAMG}.}
    In particular, in \cref{apx:spectral_equivalence} we discuss the spectral equivalence between \eqref{eq:discreteReactionDiffusion} and \eqref{eq:discreteReactionDiffusion2}.
    \bb{Throughout the finite element experiments of this section (Tables~\ref{tab:reactionDiffusionAMG}--\ref{tab:reactionDiffusionAMG3} and Figure~\ref{fig:convergenceHistory}) the quantity we monitor and on which we base the stopping criterion is the \emph{unpreconditioned} residual $\lVert A\vec{x}_k - \vec{b}\rVert_2$ of the original linear system \eqref{eq:linearsystem}, and not the preconditioned residual nor the residual of the normal equation \eqref{eq:normalPDE}. Reporting the residual of the original system makes the iteration counts directly comparable with those of the FGMRES/AIR solver applied to \eqref{eq:linearsystem} in Table~\ref{tab:air}.}
    \begin{table}
        \centering
        \caption{Comparison of the number of iterations for the CGNR method preconditioned by the inversion via PETSc GAMG of \eqref{eq:discreteReactionDiffusion2}, for different values of $\nu$ and different mesh sizes.
            The wind is fixed to $\beta = (1,0)$, and as right-hand side we consider the function $f(x,y) \equiv 1$.
            The CGNR method was terminated when the \bb{residual norm}  was less than $10^{-5}$. \bb{The discretisation is SUPG-stabilised with $\delta = h/2|\beta|$ and the surrogate \eqref{eq:discreteReactionDiffusion2} is stabilised consistently. The GAMG smoother is eight symmetric-SOR sweeps per level.}}
        \label{tab:reactionDiffusionAMG}
        \begin{tabular}{c|c|c|c|c|c}
            \toprule
            $\nu$               & $32\times 32$ & $64\times 64$ & $128\times 128$ & $256\times 256$ & $512\times 512$ \\
            \midrule
            $1\cdot 10^{-2}$    &  4 &  4 &  6 &  9 & 13 \\
            $5\cdot 10^{-3}$    &  6 &  6 &  6 &  7 & 11 \\
            $2.5\cdot 10^{-3}$  & 10 &  9 &  8 &  7 &  9 \\
            $1.25\cdot 10^{-3}$ & 14 & 13 & 12 & 10 & 10 \\
            \bottomrule
        \end{tabular}
        \vspace{1cm}
    \end{table}
    \begin{table}
        \centering
        \caption{Comparison of the number of iterations for the FGMRES method preconditioned by HYPRE, with AIR(1) restriction and one step of Jacobi relaxation up and down.
            The wind is fixed to $\beta = (1,0)$, and as right-hand side we consider the function $f(x,y) \equiv 1$.
            The FGMRES method was terminated when the absolute residual  \bb{norm} was less than $10^{-5}$, a dash denotes that the method stagned and did not \bb{converge} within 40 iterations.}
        \label{tab:air}
        \begin{tabular}{c|c|c|c|c|c}
            \toprule
            $\nu$               & $32\times 32$ & $64\times 64$ & $128\times 128$ & $256\times 256$ & $512\times 512$ \\
            \midrule
            $1\cdot 10^{-2}$    & 3             & 3             & 4               & 5               & 7               \\
            $5\cdot 10^{-3}$    & -             & 5             & 4               & 5               & 6               \\
            $2.5\cdot 10^{-3}$  & -             & -             & 6               & 5               & 5               \\
            $1.25\cdot 10^{-3}$ & -             & -             & -               & -               & 7               \\
            \bottomrule
        \end{tabular}
        \vspace{1cm}
    \end{table}
    \begin{figure}[h]
        \centering
        \includegraphics[scale=0.066]{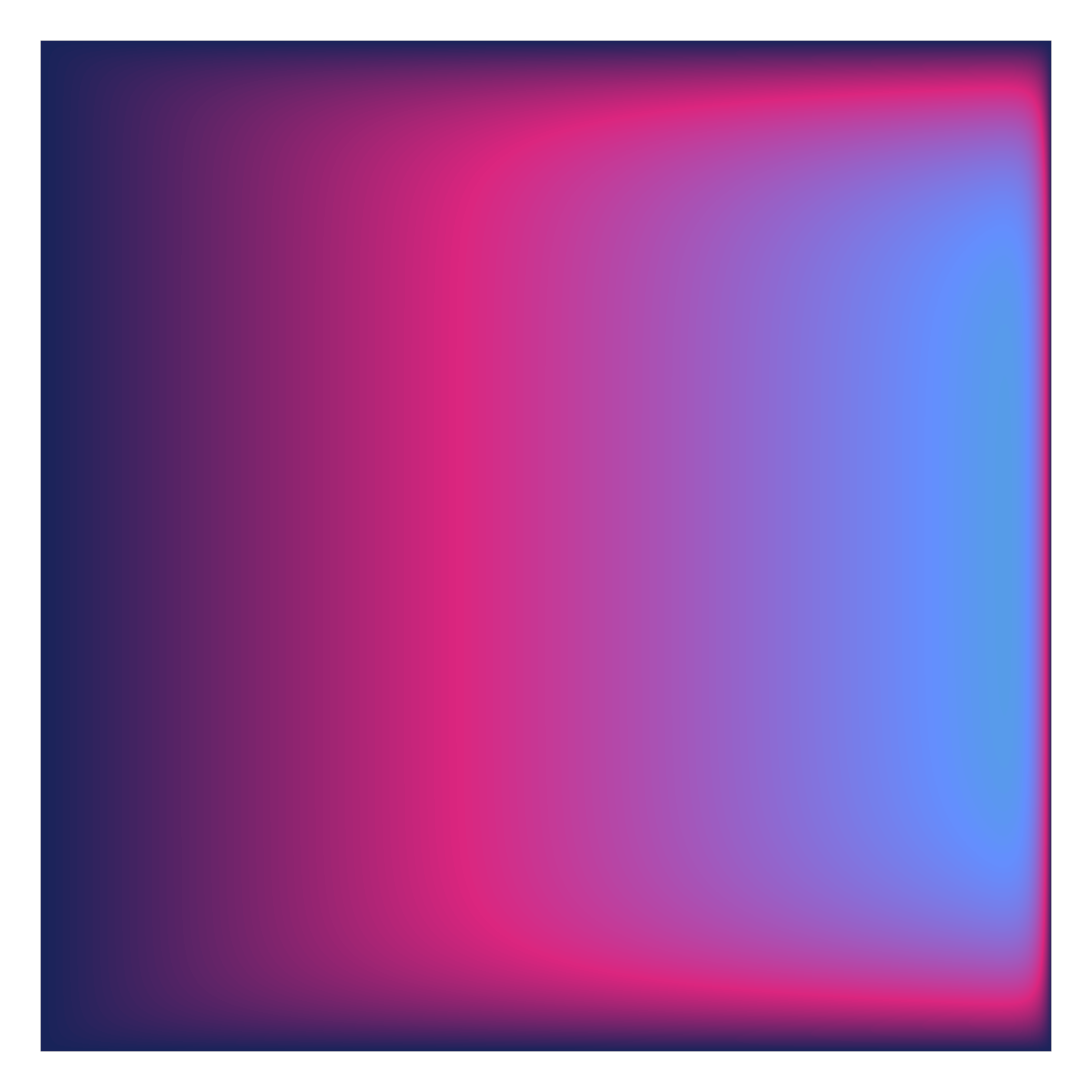}
        \includegraphics[scale=0.066]{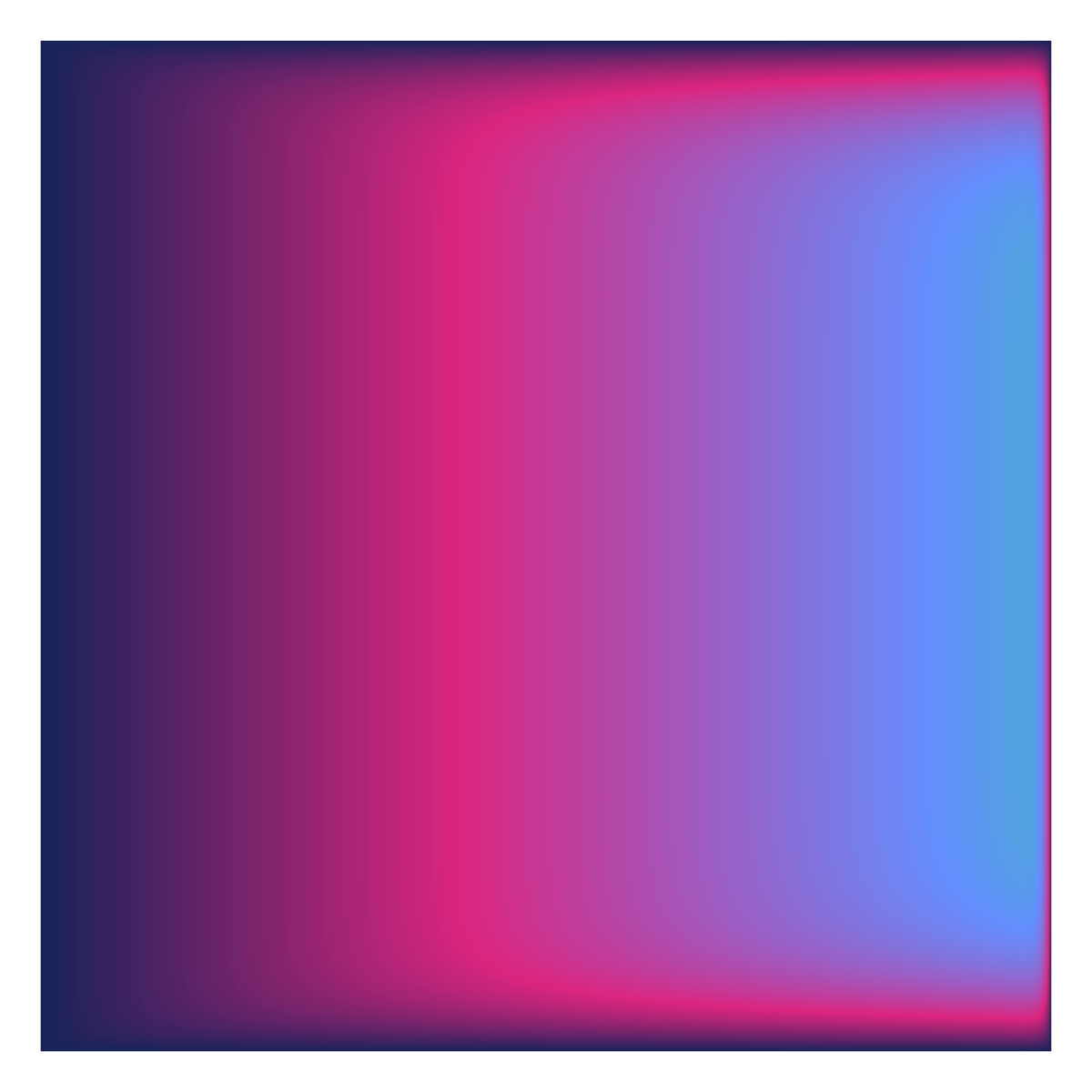}
        \includegraphics[scale=0.066]{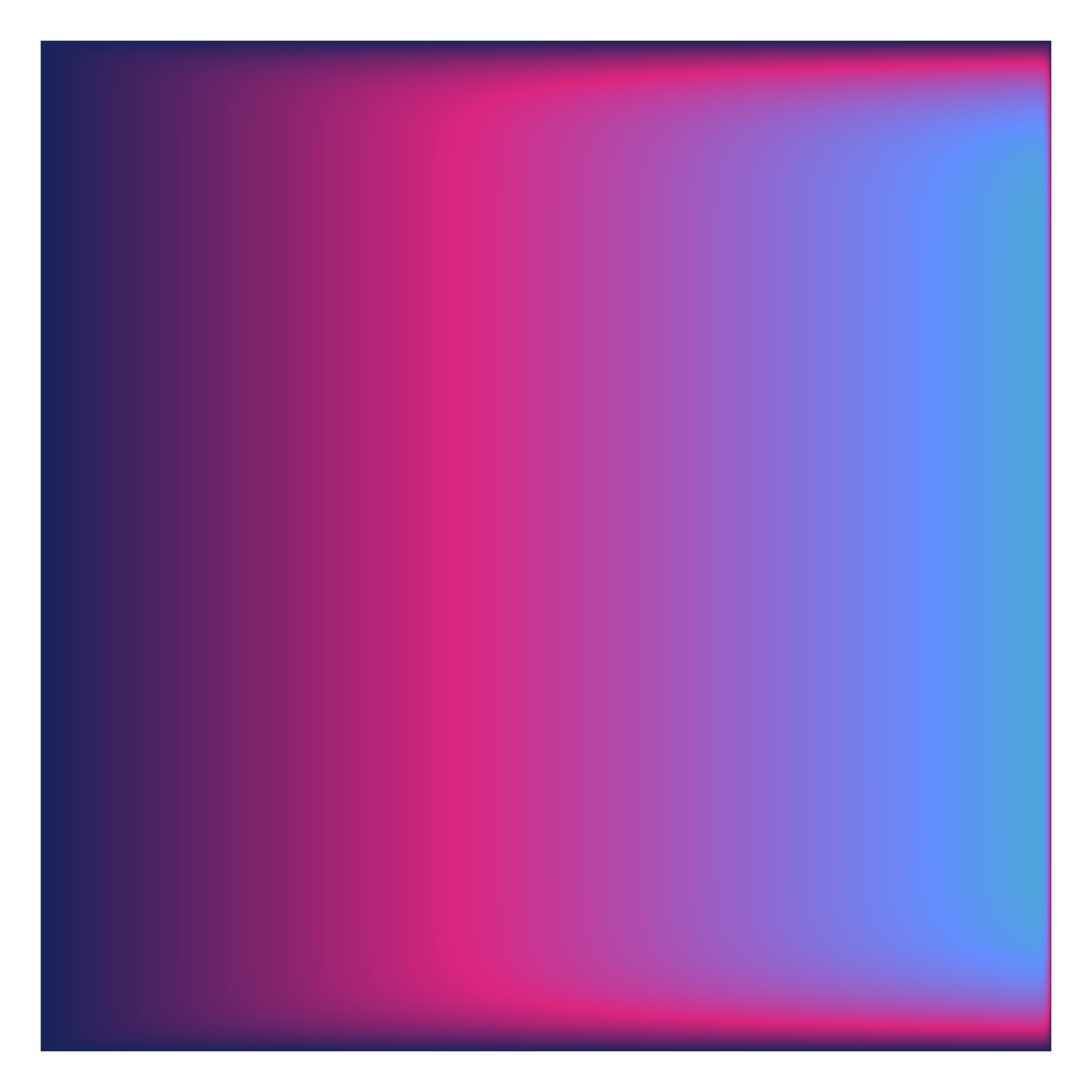}
        \includegraphics[scale=0.066]{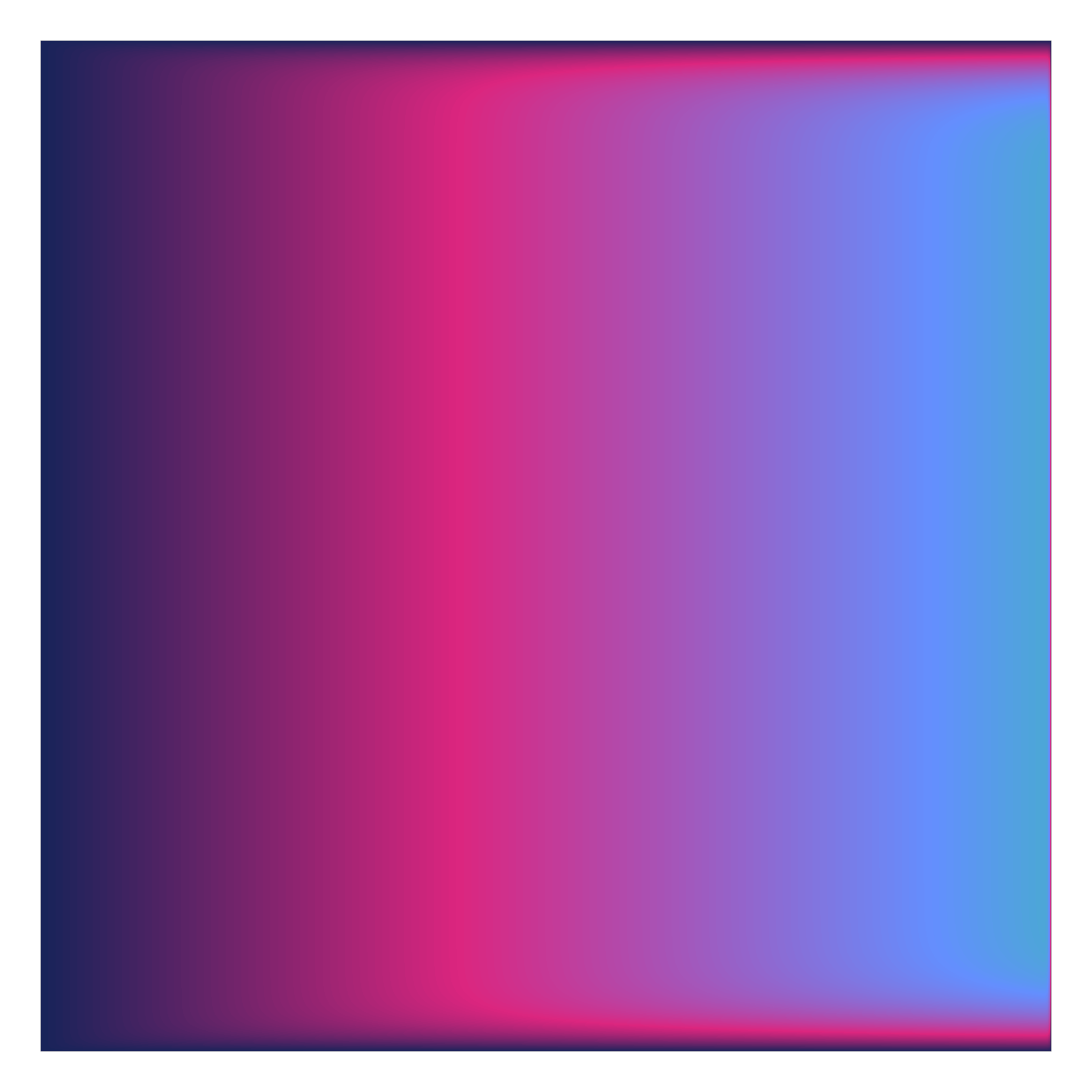}\hspace{3pt}\includegraphics[height=2.79cm]{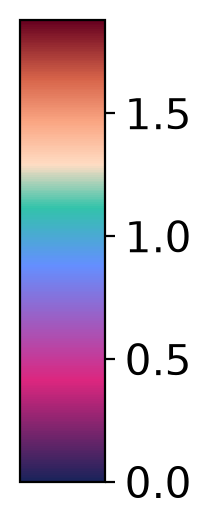}

        \vspace{0.5cm}

        \includegraphics[scale=0.066]{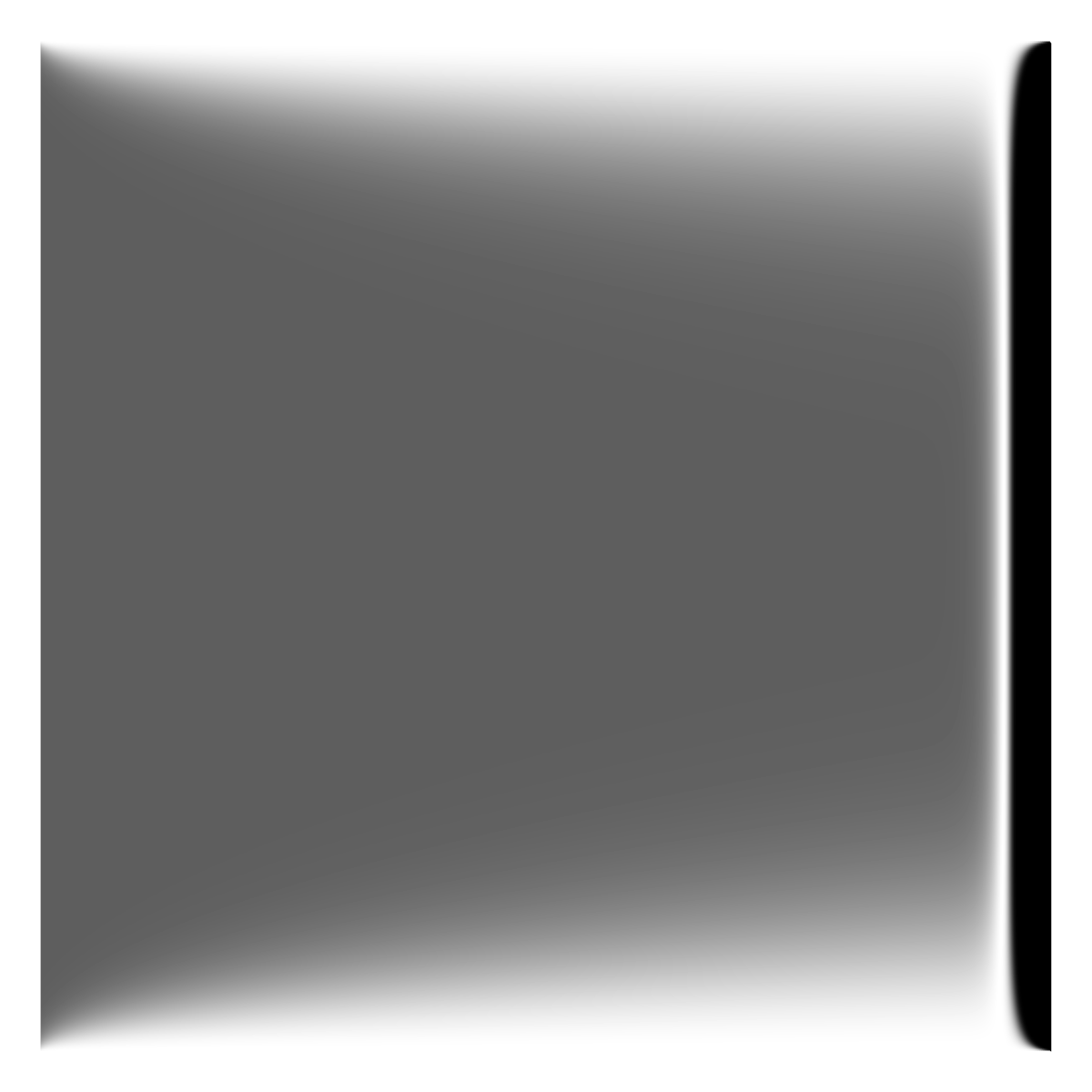}
        \includegraphics[scale=0.066]{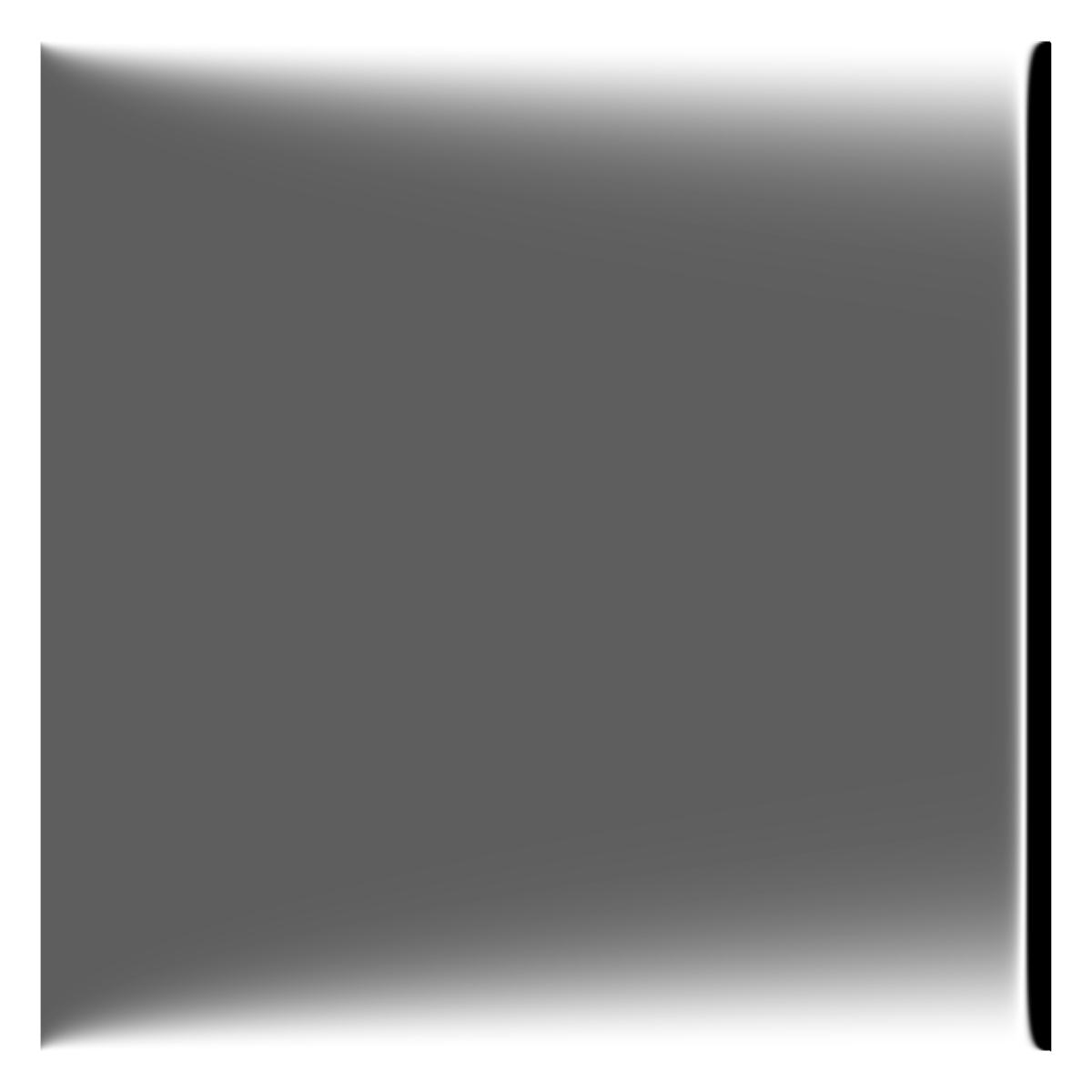}
        \includegraphics[scale=0.066]{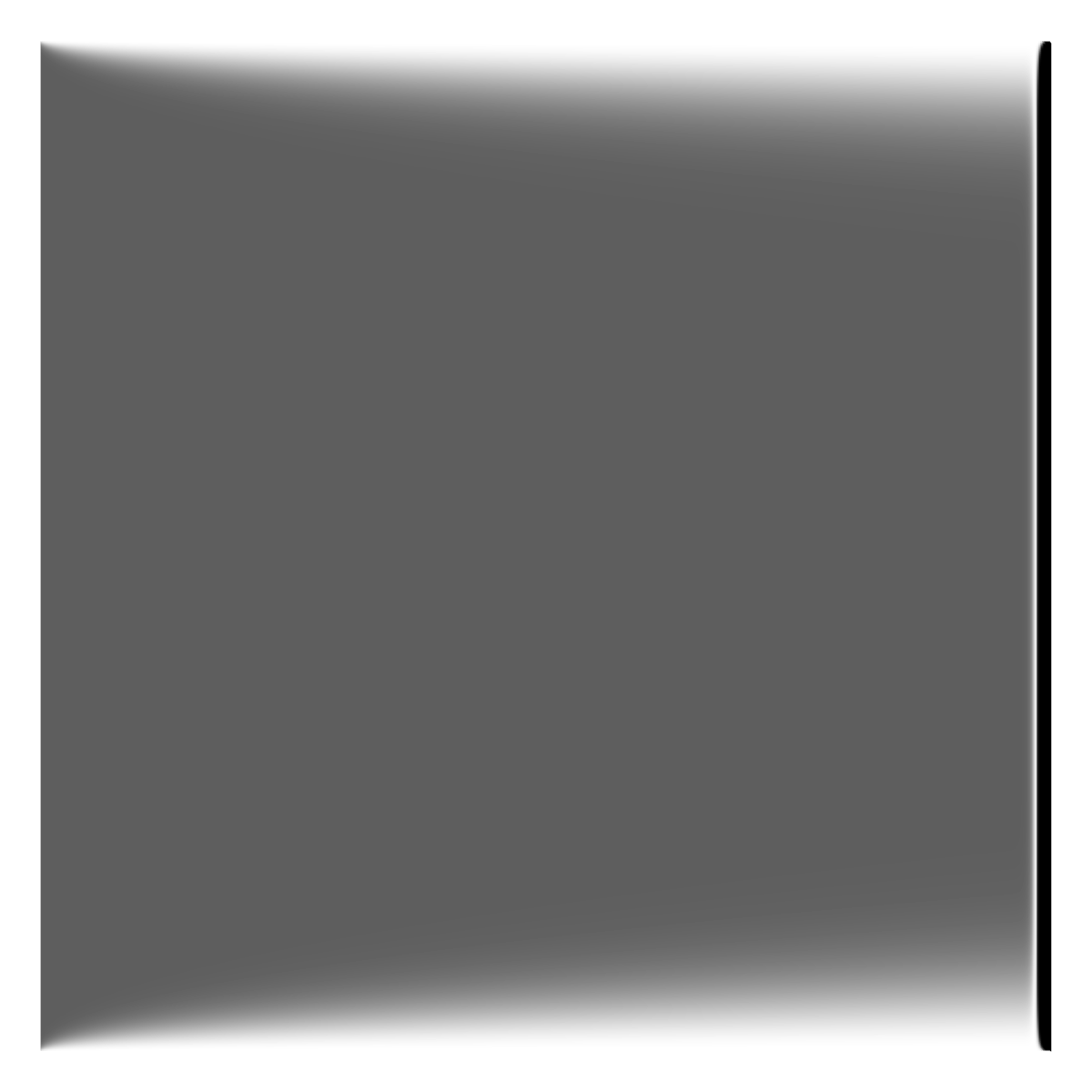}
        \includegraphics[scale=0.066]{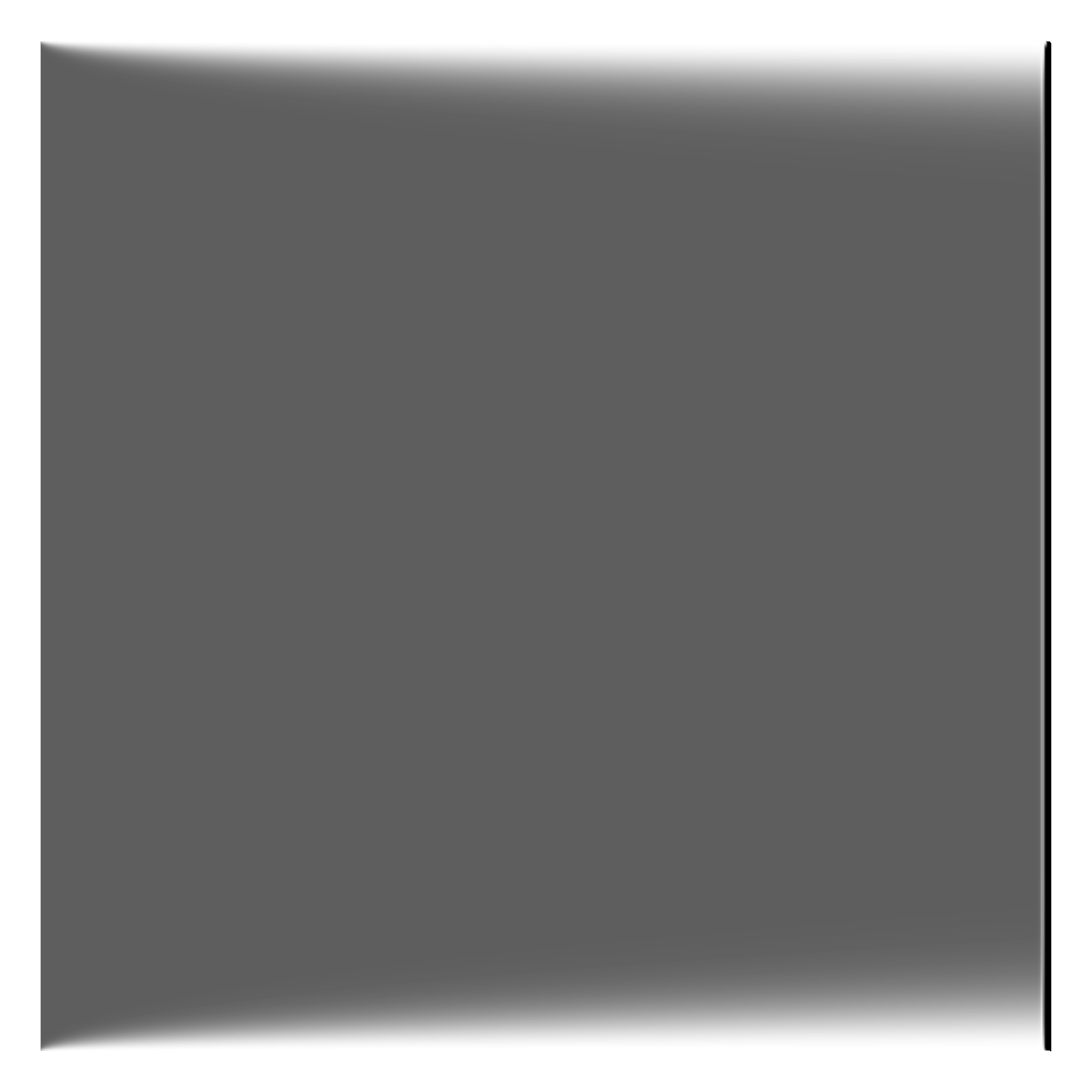}\hspace{3pt}\includegraphics[height=2.79cm]{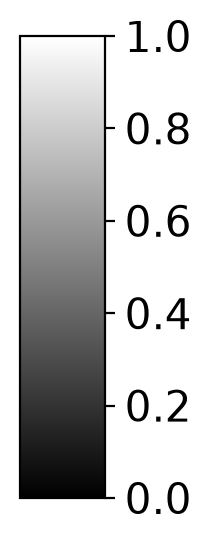}
        \caption{The discrete solution $u_h$ of the reaction diffusion equation \eqref{eq:weakForm}, with $\vec{\beta}=(1,0)$, for different values of $\nu$ at the finest mesh size $512\times 512$, together with $exp(-\lvert\nabla\cdot \beta u_h\lvert^2)$.}
        \label{fig:reactionDiffusionAMG}
    \end{figure}
    \begin{table}
        \centering
        \caption{Comparison of the number of iterations for the CGNR method preconditioned by the inversion via PETSc GAMG of \eqref{eq:discreteReactionDiffusion2}, for different values of $\nu$ and different mesh sizes.
            The wind is fixed to $\norm{\beta}\beta = (1,1)$, and as right-hand side we consider the function $f(x,y) \equiv 1$.
            The CGNR method was terminated when the absolute residual  \bb{norm} was less than $10^{-5}$. \bb{The discretisation is SUPG-stabilised with $\delta = h/2|\beta|$ and the surrogate \eqref{eq:discreteReactionDiffusion2} is stabilised consistently; the GAMG smoother is eight symmetric-SOR sweeps per level.}}
        \label{tab:reactionDiffusionAMG2}
        \begin{tabular}{c|c|c|c|c|c}
            \toprule
            $\nu$               & $32\times 32$ & $64\times 64$ & $128\times 128$ & $256\times 256$ & $512\times 512$ \\
            \midrule
            $1\cdot 10^{-2}$    &  4 &  5 &  7 & 10 & 15 \\
            $5\cdot 10^{-3}$    &  7 &  6 &  7 & 10 & 15 \\
            $2.5\cdot 10^{-3}$  & 11 & 10 &  9 & 10 & 14 \\
            $1.25\cdot 10^{-3}$ & 16 & 16 & 14 & 13 & 15 \\
            \bottomrule
        \end{tabular}
    \end{table}
    \begin{figure}[h]
        \centering
        \includegraphics[scale=0.066]{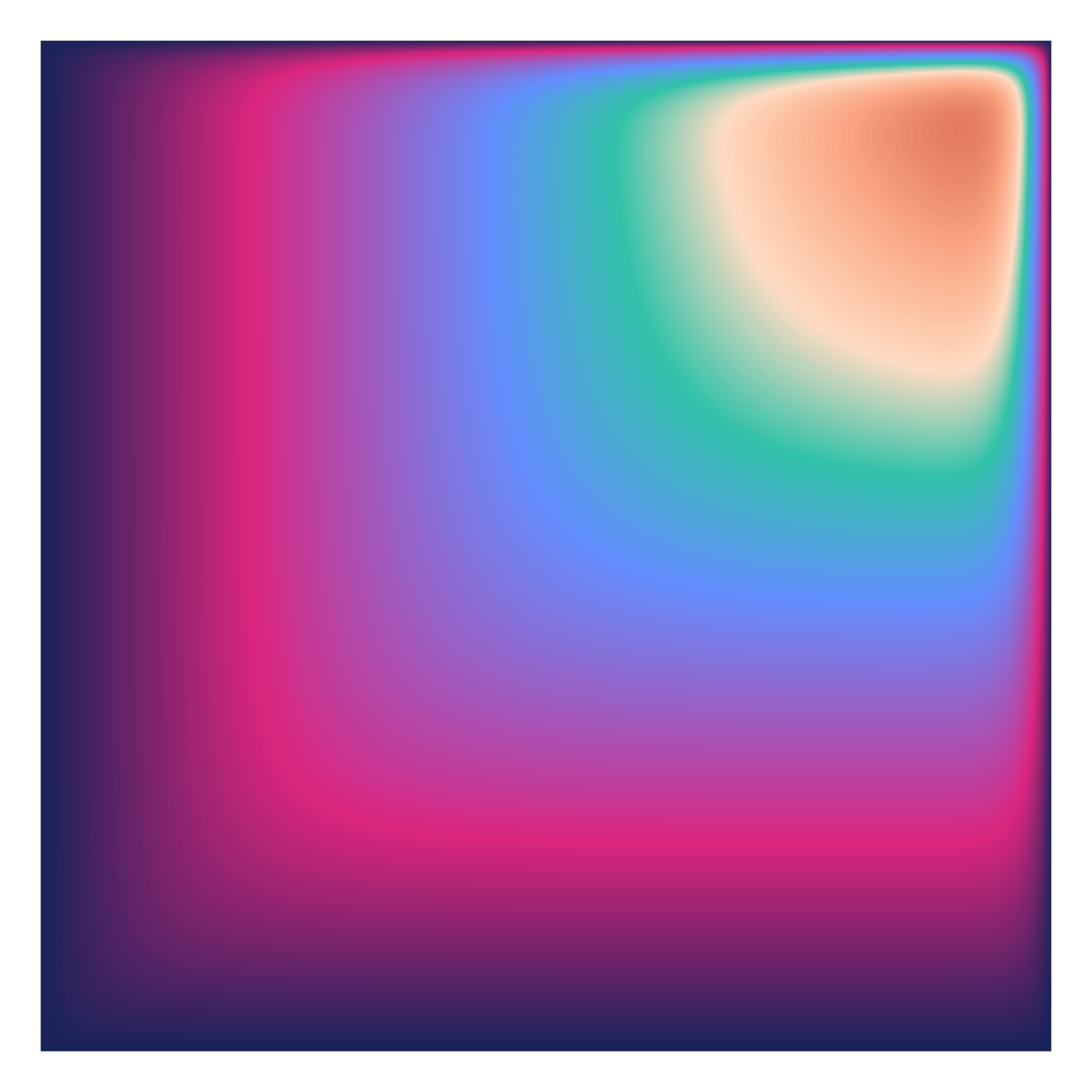}
        \includegraphics[scale=0.066]{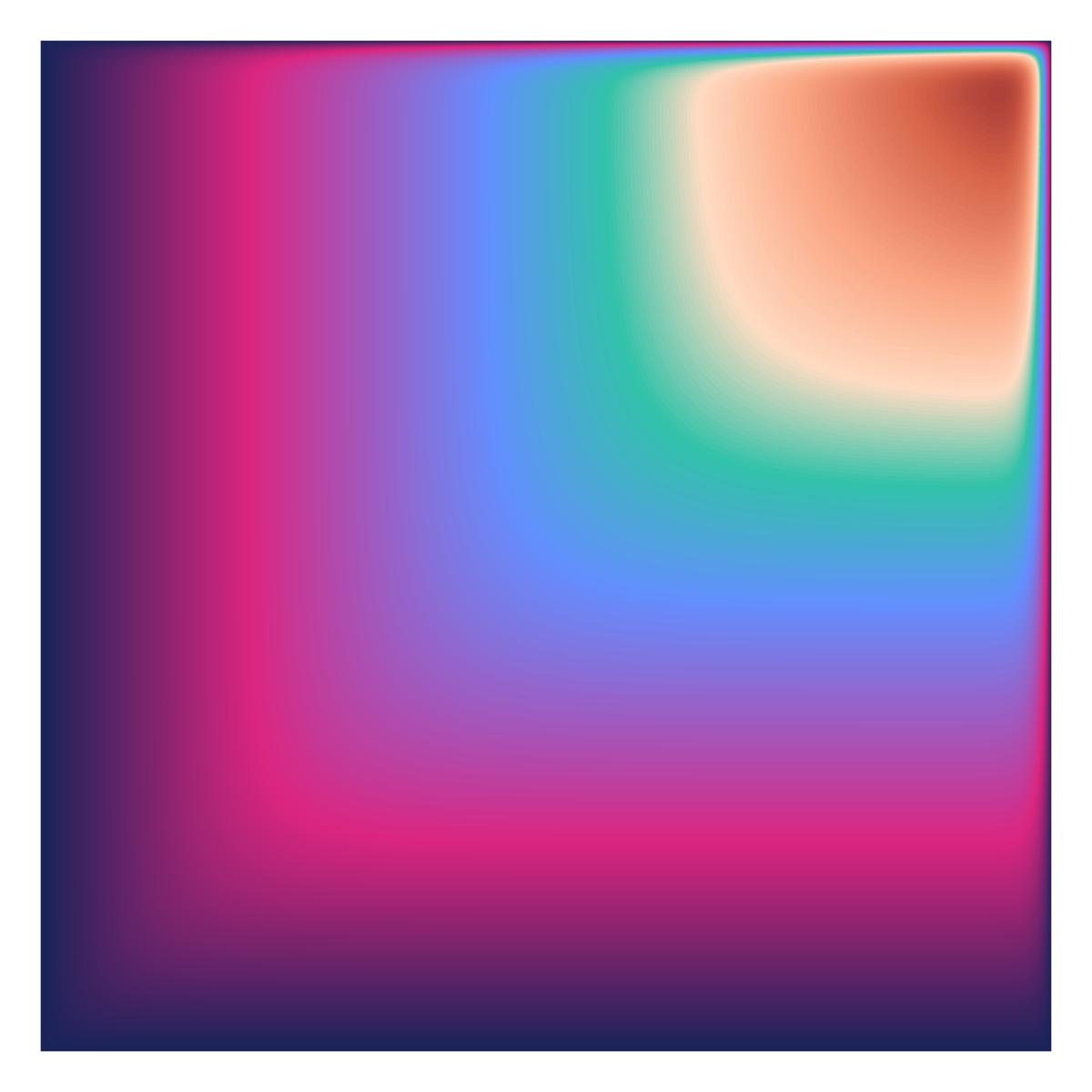}
        \includegraphics[scale=0.066]{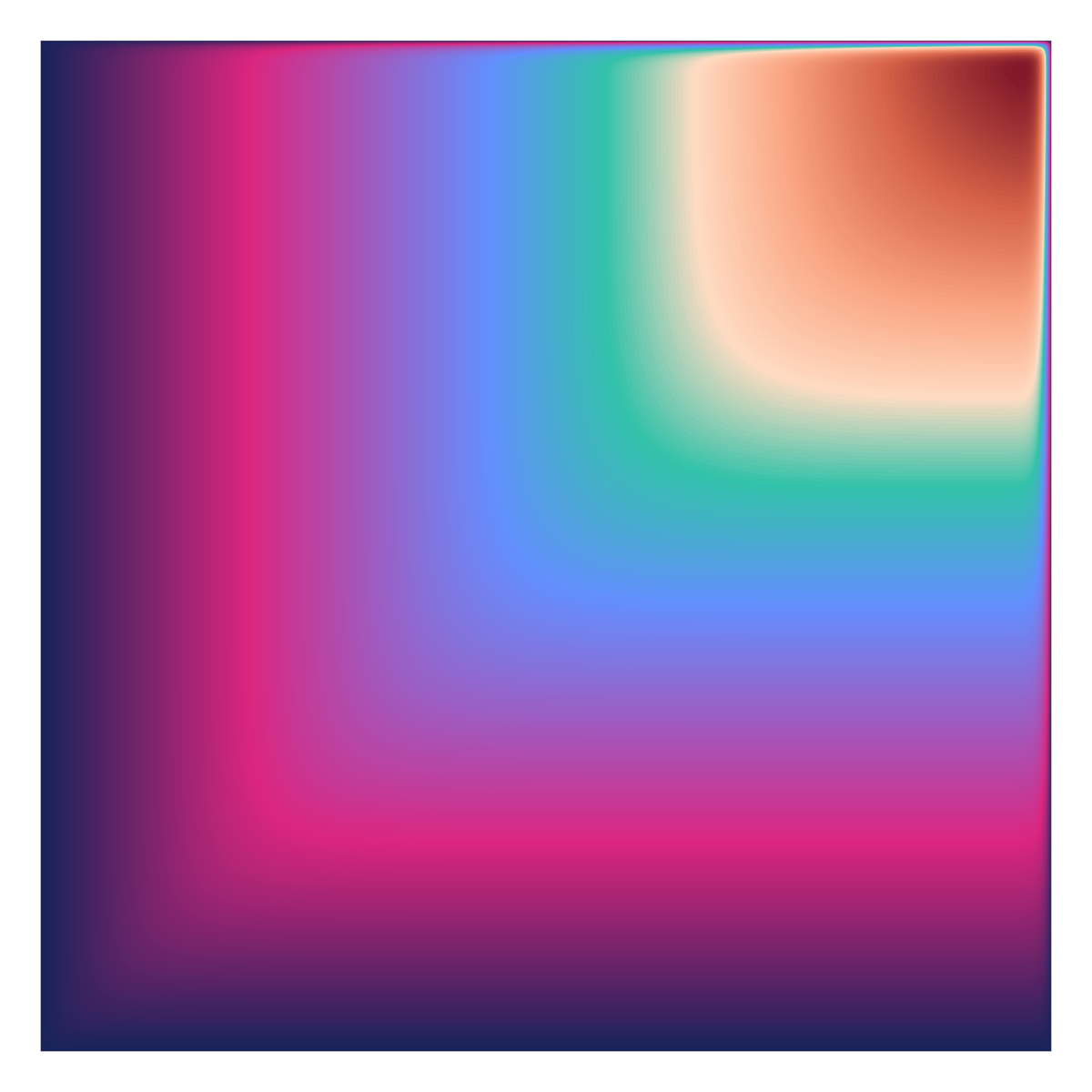}
        \includegraphics[scale=0.066]{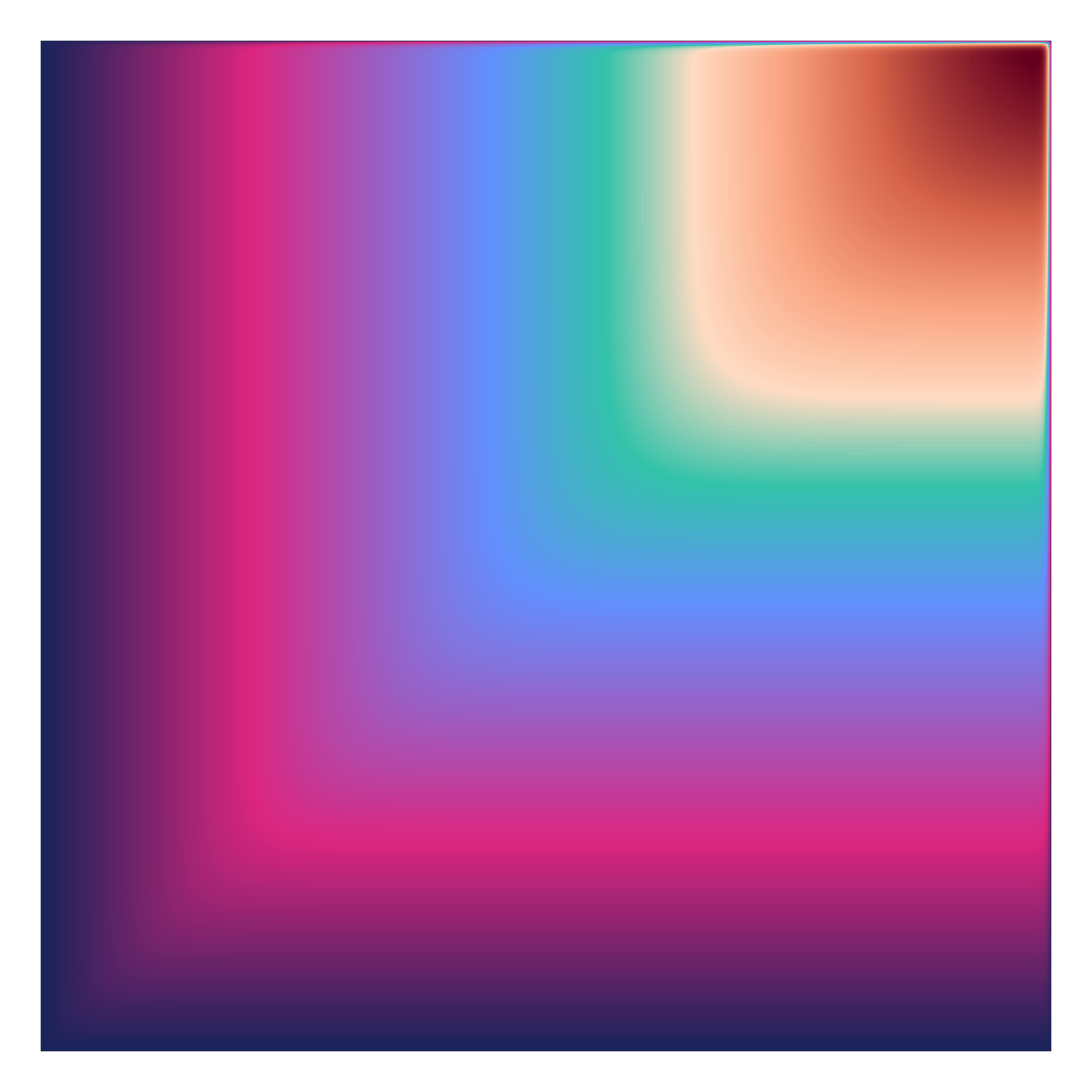}\hspace{3pt}\includegraphics[height=2.79cm]{figures/cw_sol_colorbar.png}

        \vspace{0.5cm}

        \includegraphics[scale=0.066]{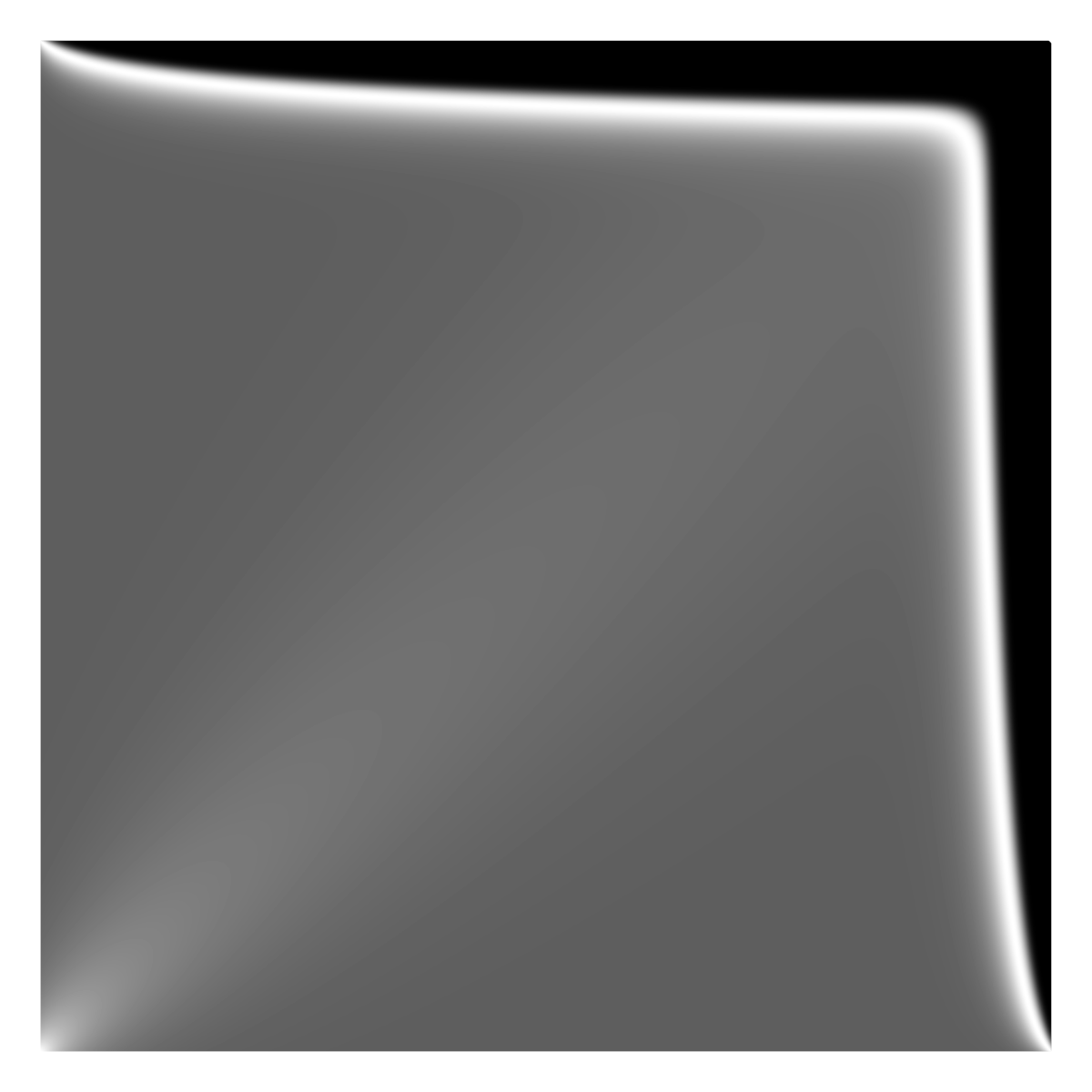}
        \includegraphics[scale=0.066]{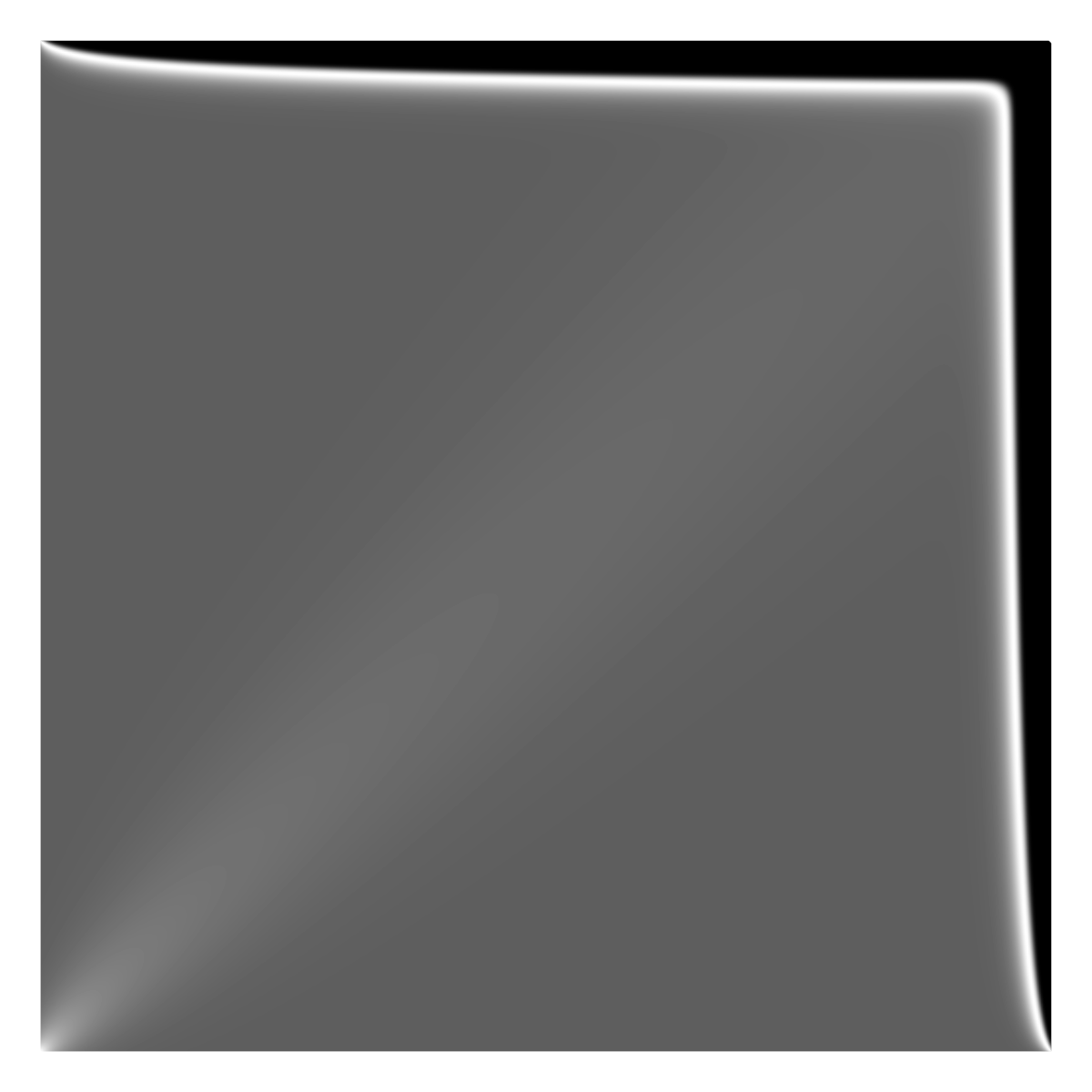}
        \includegraphics[scale=0.066]{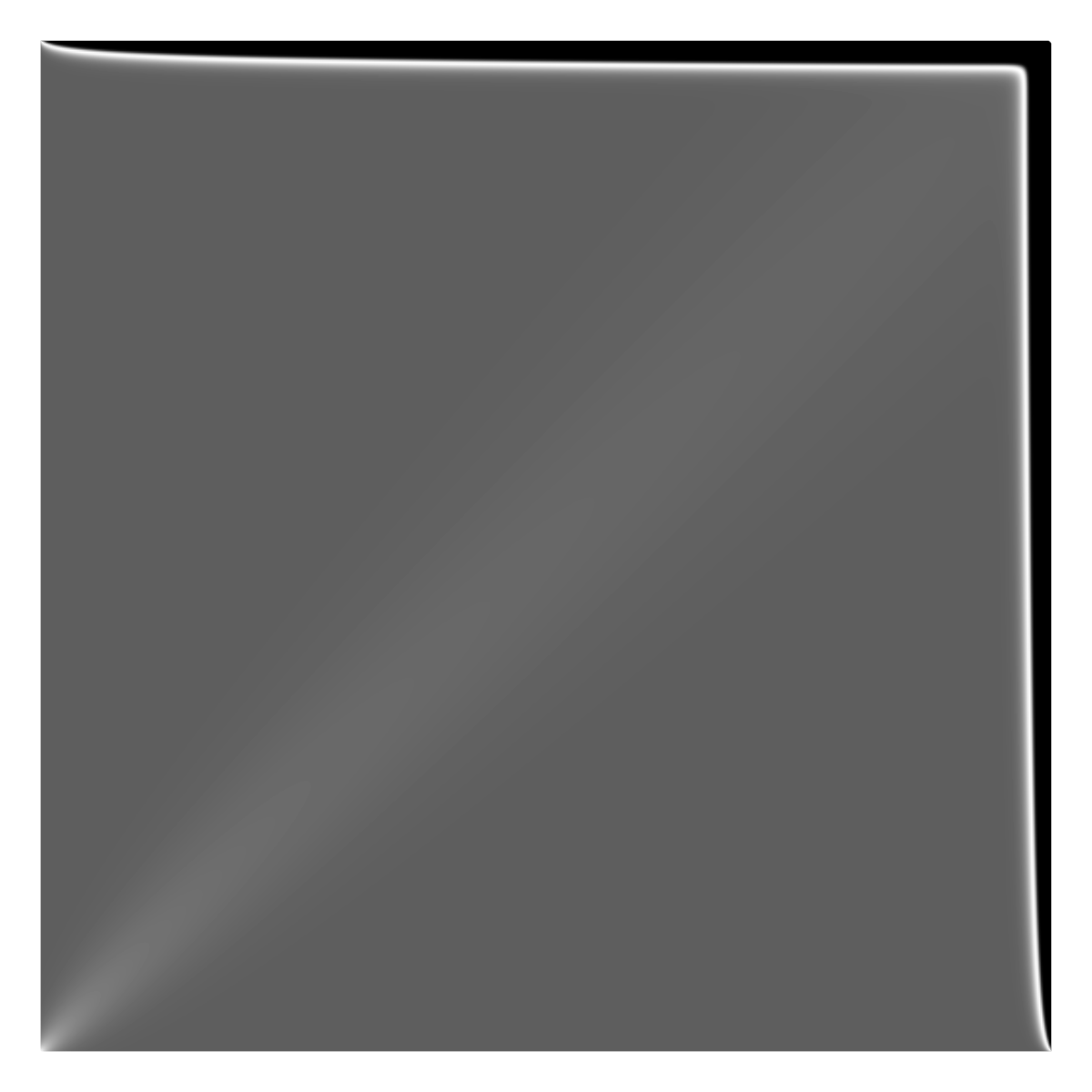}
        \includegraphics[scale=0.066]{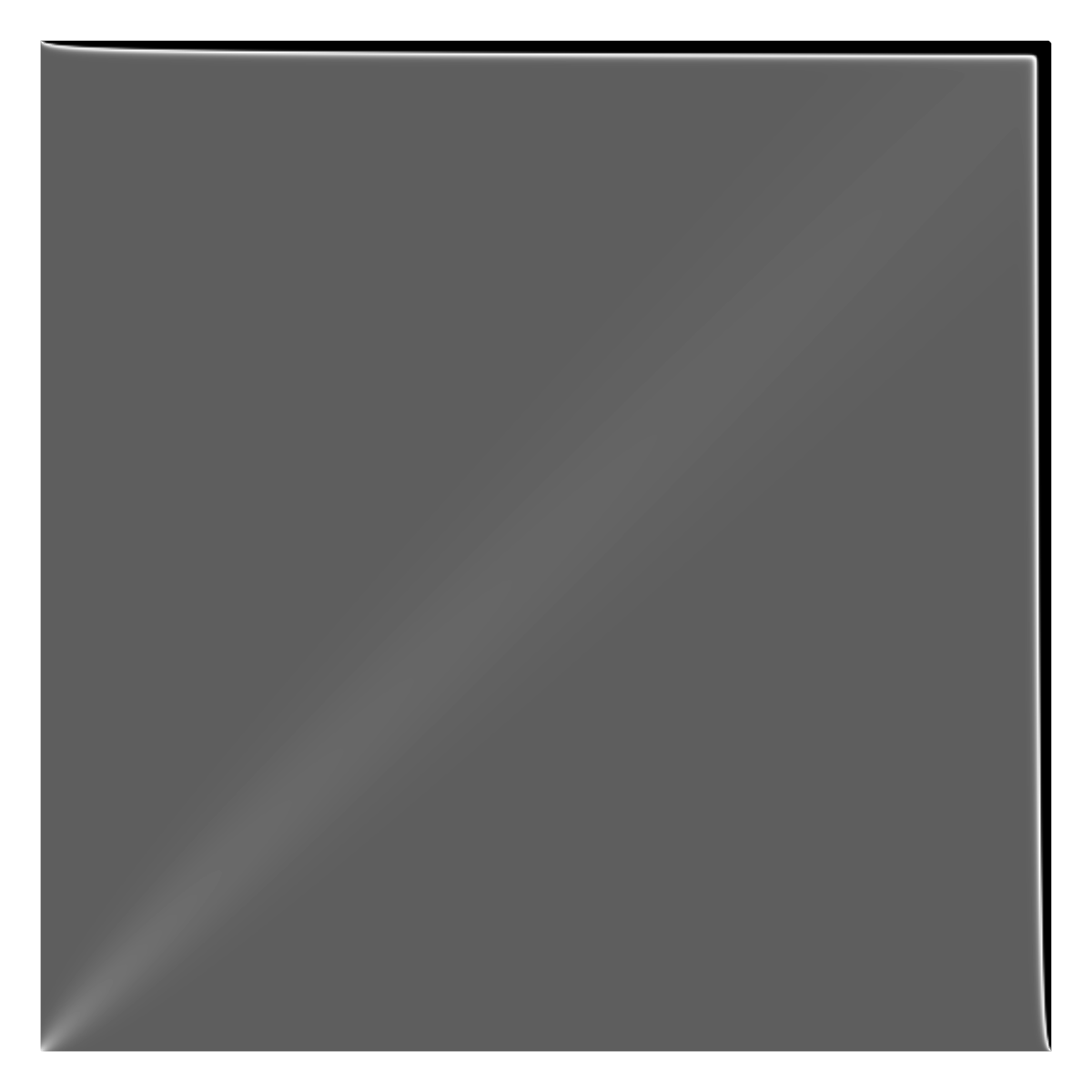}\hspace{3pt}\includegraphics[height=2.79cm]{figures/cw_div_colorbar.png}
        \caption{The discrete solution $u_h$ of the reaction diffusion equation \eqref{eq:weakForm}, with $\sqrt{2}\vec{\beta}=(1,1)$, for different values of $\nu$ at the finest mesh size $512\times 512$, together with $exp(-\lvert\nabla\cdot \beta u_h\lvert^2)$.}
        \label{fig:reactionDiffusionAMG2}
    \end{figure}
    \begin{figure}[h]
        \centering
        \begin{tikzpicture}[scale=0.75]
            \begin{groupplot}[%
                    group style={%
                            group size=2 by 2,
                            horizontal sep=1.5cm,
                            vertical sep=2cm,
                        },
                    ymajorgrids=true,
                    grid style=dashed,
                ]
                \nextgroupplot[width=8cm,height=6cm,domain=2:5,ymode=log, xlabel={Krylov steps}, ylabel={$\norm{Ax-b}_{2}$}, title={$32\times 32$}, 
                    legend style={legend columns=3, draw=none,nodes={scale=.8}},
                    legend to name=name]
                \addplot+[line width=1.5pt,mark=None] table [x=step, y=res, col sep=comma] {data/fem_advection_mass_normal_eq_diagflow_0.00125_32.0.csv};
                \addplot+[line width=1.5pt,mark=None, dashed] table [x=step, y=res, col sep=comma] {data/fem_advection_mass_normal_eq_diagflow_0.0025_32.0.csv};
                \addplot+[line width=1.5pt,mark=None, dotted] table [x=step, y=res, col sep=comma] {data/fem_advection_mass_normal_eq_diagflow_0.01_32.0.csv};
                \nextgroupplot[width=8cm,height=6cm,domain=2:5,ymode=log, xlabel={Krylov steps}, ylabel={$\,$}, title={$64\times 64$}, 
                    legend style={legend columns=3, draw=none,nodes={scale=.8}},
                    legend to name=named]
                \addplot+[line width=1.5pt,mark=None] table [x=step, y=res, col sep=comma] {data/fem_advection_mass_normal_eq_diagflow_0.00125_64.0.csv};
                \addplot+[line width=1.5pt,mark=None, dashed] table [x=step, y=res, col sep=comma] {data/fem_advection_mass_normal_eq_diagflow_0.0025_64.0.csv};
                \addplot+[line width=1.5pt,mark=None, dotted] table [x=step, y=res, col sep=comma] {data/fem_advection_mass_normal_eq_diagflow_0.01_64.0.csv};
                \nextgroupplot[width=8cm,height=6cm,domain=2:5,ymode=log, xlabel={Krylov steps}, ylabel={$\norm{Ax-b}_{2}$}, title={$128\times 128$}, 
                    legend style={legend columns=3, draw=none,nodes={scale=.8}},
                    legend to name=named]
                \addplot+[line width=1.5pt,mark=None] table [x=step, y=res, col sep=comma] {data/fem_advection_mass_normal_eq_diagflow_0.00125_128.0.csv};
                \addplot+[line width=1.5pt,mark=None, dashed] table [x=step, y=res, col sep=comma] {data/fem_advection_mass_normal_eq_diagflow_0.0025_128.0.csv};
                \addplot+[line width=1.5pt,mark=None, dotted] table [x=step, y=res, col sep=comma] {data/fem_advection_mass_normal_eq_diagflow_0.01_128.0.csv};
                \nextgroupplot[width=8cm,height=6cm,domain=2:5,ymode=log, xlabel={Krylov steps}, ylabel={$\,$}, title={$256\times 256$}, 
                    legend style={legend columns=3, draw=none,nodes={scale=.8}},
                    legend to name=named]
                \addplot+[line width=1.5pt,mark=None] table [x=step, y=res, col sep=comma] {data/fem_advection_mass_normal_eq_diagflow_0.00125_256.0.csv};
                \addplot+[line width=1.5pt,mark=None, dashed] table [x=step, y=res, col sep=comma] {data/fem_advection_mass_normal_eq_diagflow_0.0025_256.0.csv};
                \addplot+[line width=1.5pt,mark=None, dotted] table [x=step, y=res, col sep=comma] {data/fem_advection_mass_normal_eq_diagflow_0.01_256.0.csv};
                \legend{$\nu = 1.25 \cdot 10^{-3}$,$\nu = 2.5 \cdot 10^{-3}$,$\nu = 1 \cdot 10^{-2}$}
            \end{groupplot}
        \end{tikzpicture}
        \pgfplotslegendfromname{named}
        \caption{Comparison of convergence history for the CGNR method preconditioned by the inversion via PETSc GAMG of \eqref{eq:discreteReactionDiffusion2}, for different values of $\nu$ and different mesh sizes.
            The wind is fixed to $\norm{\beta}\beta = (1,1)$, and as right-hand side we consider the function $f(x,y) \equiv 1$.
            The CGNR method was terminated when the absolute residual  \bb{norm} was less than $10^{-5}$.}
        \label{fig:convergenceHistory}
    \end{figure}
    Observing Table \ref{tab:reactionDiffusionAMG} and Table \ref{tab:reactionDiffusionAMG2}, we see that the approach here presented yields an effective preconditioner for different types of constant winds, i.e.~ $\beta=(1,0)$ and $\beta=(1,1)$, respectively depicted in \cref{fig:reactionDiffusionAMG} and \cref{fig:reactionDiffusionAMG2}.
    {Throughout these experiments the preconditioner is a single GAMG V-cycle on the sparse unprojected operator \eqref{eq:discreteReactionDiffusion2}, stabilised consistently with the SUPG term carried by $A$, with the smoother tuned to eight symmetric-SOR sweeps per level.}
    In particular, the convergence history for iteration counts presented in Table \ref{tab:reactionDiffusionAMG2}, is displayed in Figure \ref{fig:convergenceHistory}.
    This can be compared to Table \ref{tab:air} where we use FGMRES preconditioned via the algebraic multigrid preconditioner HYPRE with Approximate Ideal Restriction (AIR) \cite{air,air2,air3} directly on \eqref{eq:linearsystem}. The proposed preconditioning strategy consistently reaches the desired accuracy within 40 iterations, whereas FGMRES preconditioned via HYPRE with AIR relaxation, although converging in very few iterations on a sufficiently fine mesh, often fails to achieve the desired accuracy within 40 iterations, in addition to relying on a preconditioner technically more complicated than a classical smoothed aggregation algebraic multigrid.

\bb{To address whether the discrete solution resolves the boundary and internal layers as $\nu$ decreases and coarse meshes are used, Figure~\ref{fig:layers} displays the discrete solution $u_h$ of the $\beta=(1,0)$ problem on every mesh from $32\times 32$ to $256\times 256$ and for every viscosity from $10^{-2}$ down to $1.25\cdot 10^{-3}$. With the SUPG stabilisation $\delta = h/2|\beta|$ the outflow layer, whose physical width scales like $\nu$, is rendered \emph{monotone and oscillation-free} on every mesh: on the coarsest meshes it is spread over a few cells by the mesh-proportional streamline diffusion, and since that added diffusion scales like $h$ it recedes and the layer sharpens to its physical width as the mesh is refined. We emphasise that the mesh-independence of the iteration counts reported above does \emph{not} require the layer to be fully resolved: because the surrogate \eqref{eq:discreteReactionDiffusion2} is stabilised consistently with $A$, it stays spectrally equivalent to the stabilised normal-PDE operator on the coarse, under-resolved meshes too, i.e.~exactly the regime in which convection--diffusion preconditioners are usually most fragile.}

\newcommand{\layerpanel}[2]{%
    \raisebox{-0.5\height}{\IfFileExists{figures/layers_nu#1_n#2.png}%
        {\includegraphics[width=3cm]{figures/layers_nu#1_n#2.png}}%
        {\fbox{\parbox[c][2.5cm][c]{2.8cm}{\centering\scriptsize $u_h$\\[2pt]$\nu{=}#1$,\ $#2{\times}#2$}}}}%
}
\begin{figure}[htbp]
    \centering
    \setlength{\tabcolsep}{2pt}
    \renewcommand{\arraystretch}{1.1}
    \resizebox{\textwidth}{!}{{%
    \begin{tabular}{ccccc@{\hskip 6pt}c}
                    & $32\times 32$ & $64\times 64$ & $128\times 128$ & $256\times 256$ & \\
        \raisebox{-0.5\height}{$\nu=10^{-2}$}        & \layerpanel{1e-2}{32}    & \layerpanel{1e-2}{64}    & \layerpanel{1e-2}{128}    & \layerpanel{1e-2}{256}    & \raisebox{-0.5\height}{\includegraphics[height=3cm]{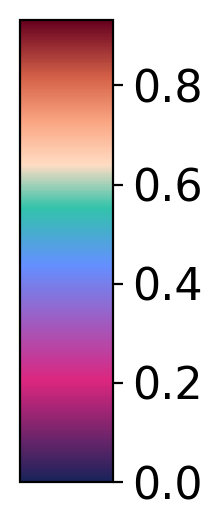}} \\
        \raisebox{-0.5\height}{$\nu=5\cdot 10^{-3}$} & \layerpanel{5e-3}{32}    & \layerpanel{5e-3}{64}    & \layerpanel{5e-3}{128}    & \layerpanel{5e-3}{256}    & \raisebox{-0.5\height}{\includegraphics[height=3cm]{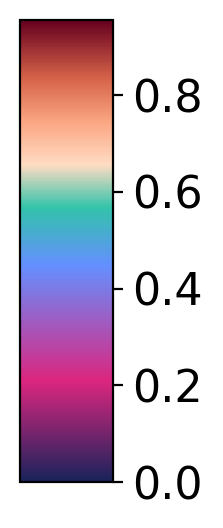}} \\
        \raisebox{-0.5\height}{$\nu=2.5\cdot 10^{-3}$} & \layerpanel{2.5e-3}{32} & \layerpanel{2.5e-3}{64} & \layerpanel{2.5e-3}{128} & \layerpanel{2.5e-3}{256} & \raisebox{-0.5\height}{\includegraphics[height=3cm]{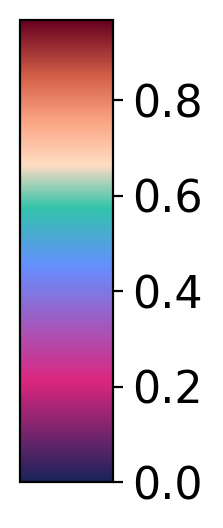}} \\
        \raisebox{-0.5\height}{$\nu=1.25\cdot 10^{-3}$} & \layerpanel{1.25e-3}{32} & \layerpanel{1.25e-3}{64} & \layerpanel{1.25e-3}{128} & \layerpanel{1.25e-3}{256} & \raisebox{-0.5\height}{\includegraphics[height=3cm]{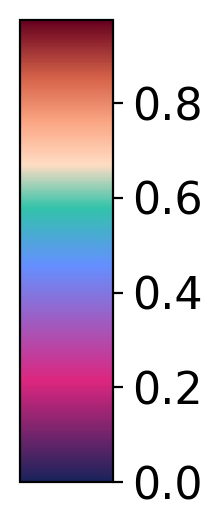}} \\
    \end{tabular}}}
    \caption{\bb{The discrete solution $u_h$ of the advection diffusion problem \eqref{eq:weakForm} with wind $\beta=(1,0)$ and $f\equiv 1$, SUPG-stabilised with $\delta = h/2|\beta|$, for each viscosity $\nu$ (rows) and each mesh size (columns). The outflow layer is monotone and oscillation-free on every mesh; on the coarse meshes it is spread over a few cells by the mesh-proportional streamline diffusion and sharpens to its physical width under refinement, whereas the iteration counts of Table~\ref{tab:reactionDiffusionAMG} remain mesh independent.}}
    \label{fig:layers}
\end{figure}

\begin{remark}\label{rmk:boundaryLayers}
    \bb{The differences in iteration count between the various configurations can be traced back to the boundary layer structure of the problem. The preconditioner \eqref{eq:discreteReactionDiffusion2} is obtained from the exact normal PDE \eqref{eq:discreteReactionDiffusion} by replacing the projected reaction term $\nu^{-1}(\Pi_{\nabla}(\beta u_h),\Pi_{\nabla}(\beta v_h))$ with the unprojected one $\nu^{-1}(\beta u_h,\beta v_h)$. The two coincide wherever $\beta u_h$ is already close to a gradient of an $H^1_0(\Omega)$ function, i.e.~almost everywhere \emph{outside} the boundary and internal layers, and differ appreciably only \emph{inside} the layers, where $\Pi_{\nabla}$ removes the layer-induced solenoidal component of $\beta u_h$. Consequently, the thinner the layers, the more accurate the approximation \eqref{eq:discreteReactionDiffusion2} of \eqref{eq:discreteReactionDiffusion}, and the tighter the spectral equivalence $AP^{-1}$ between the two operators. This explains, for instance, why imposing homogeneous Neumann conditions on the top and bottom boundaries (Table~\ref{tab:reactionDiffusionAMG3}), which suppresses the Dirichlet layers present in the pure-Dirichlet problem (Table~\ref{tab:reactionDiffusionAMG}), systematically lowers the iteration count and, more generally, why wind configurations that produce thinner layers are preconditioned more effectively by \eqref{eq:discreteReactionDiffusion2}.}
\end{remark}

    To conclude this example we consider the advection diffusion problem with mixed boundary conditions, i.e. we impose homogeneous Dirichlet boundary conditions on the left and right boundary of the unit square, while we impose homogeneous Neumann boundary conditions on the top and bottom boundary.
    Once again, we adopted a $\mathcal{P}^1$ finite element discretisation on a uniform triangulation, obtained by splitting each square of a uniform grid into two triangles cutting along the south-west to north-east diagonal.
    We inverted the discretisation of the normal PDE \eqref{eq:discreteReactionDiffusion2} with mixed boundary condition, using as a preconditioner via PETSc GAMG \cite{Adams}, a smoothed aggregation algebraic multigrid implementation. The result are displayed in Table \ref{tab:reactionDiffusionAMG3}.
    \bb{For the wind $\beta=(1,0)$ the original problem has its inflow on the left boundary and its outflow on the right, rhus the adjoint problem, governed by $-\beta\cdot\nabla$, reverses the wind, so its inflow is the outflow of the original problem and vice versa. A boundary layer of one problem at the outflow is therefore matched by a boundary layer of the other at the same location, and to resolve both operators any mesh refinement must be placed at \emph{both} the left and the right boundaries. The wind is parallel to the top and bottom boundaries so no advective boundary layer forms there for either $A$ or $A^T$.}
    \begin{table}
        \centering
        \caption{Comparison of the number of iterations for the CGNR method preconditioned by the inversion via PETSc GAMG of the reaction diffusion problem with mixed boundary conditions, for different values of $\nu$ and different mesh sizes.
            The wind is fixed to $\beta = (1,0)$, and as right-hand side we consider the function $f(x,y) \equiv 1$.
            The CGNR method was terminated when the absolute unpreconditioned residual  \bb{norm} was less than $10^{-5}$. \bb{The discretisation is SUPG-stabilised with $\delta = h/2|\beta|$ and the surrogate \eqref{eq:discreteReactionDiffusion2} is stabilised consistently; the GAMG smoother is eight symmetric-SOR sweeps per level.}}
        \label{tab:reactionDiffusionAMG3}
        \begin{tabular}{c|c|c|c|c|c}
            \toprule
            $\nu$               & $32\times 32$ & $64\times 64$ & $128\times 128$ & $256\times 256$ & $512\times 512$ \\
            \midrule
            $1\cdot 10^{-2}$    &  4 &  4 &  5 &  6 &  6 \\
            $5\cdot 10^{-3}$    &  6 &  5 &  5 &  5 &  6 \\
            $2.5\cdot 10^{-3}$  &  9 &  7 &  5 &  5 &  5 \\
            $1.25\cdot 10^{-3}$ & 11 &  9 &  7 &  5 &  5 \\
            \bottomrule
        \end{tabular}
    \end{table}
    \begin{figure}[h]
        \centering
        \includegraphics[scale=0.066]{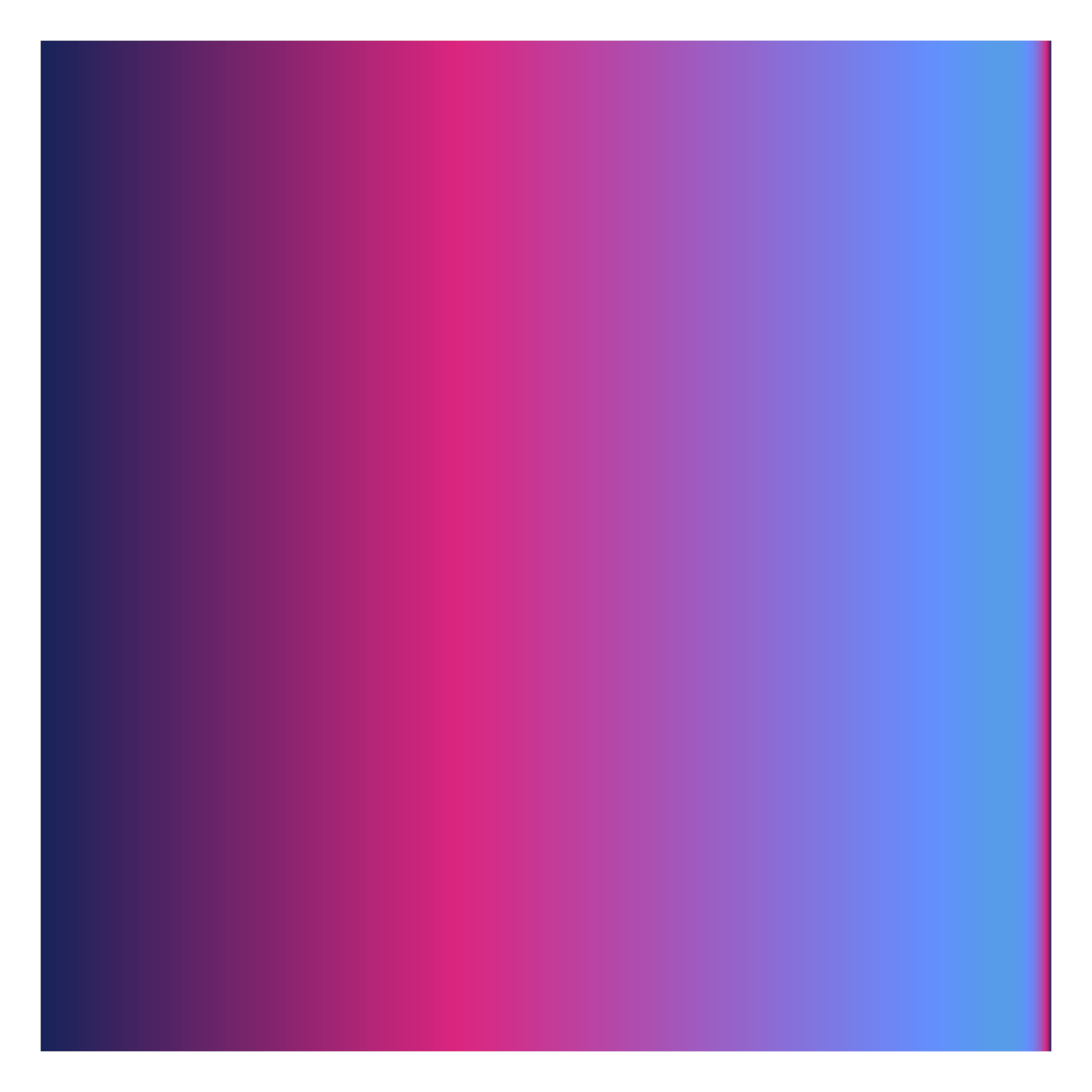}
        \includegraphics[scale=0.066]{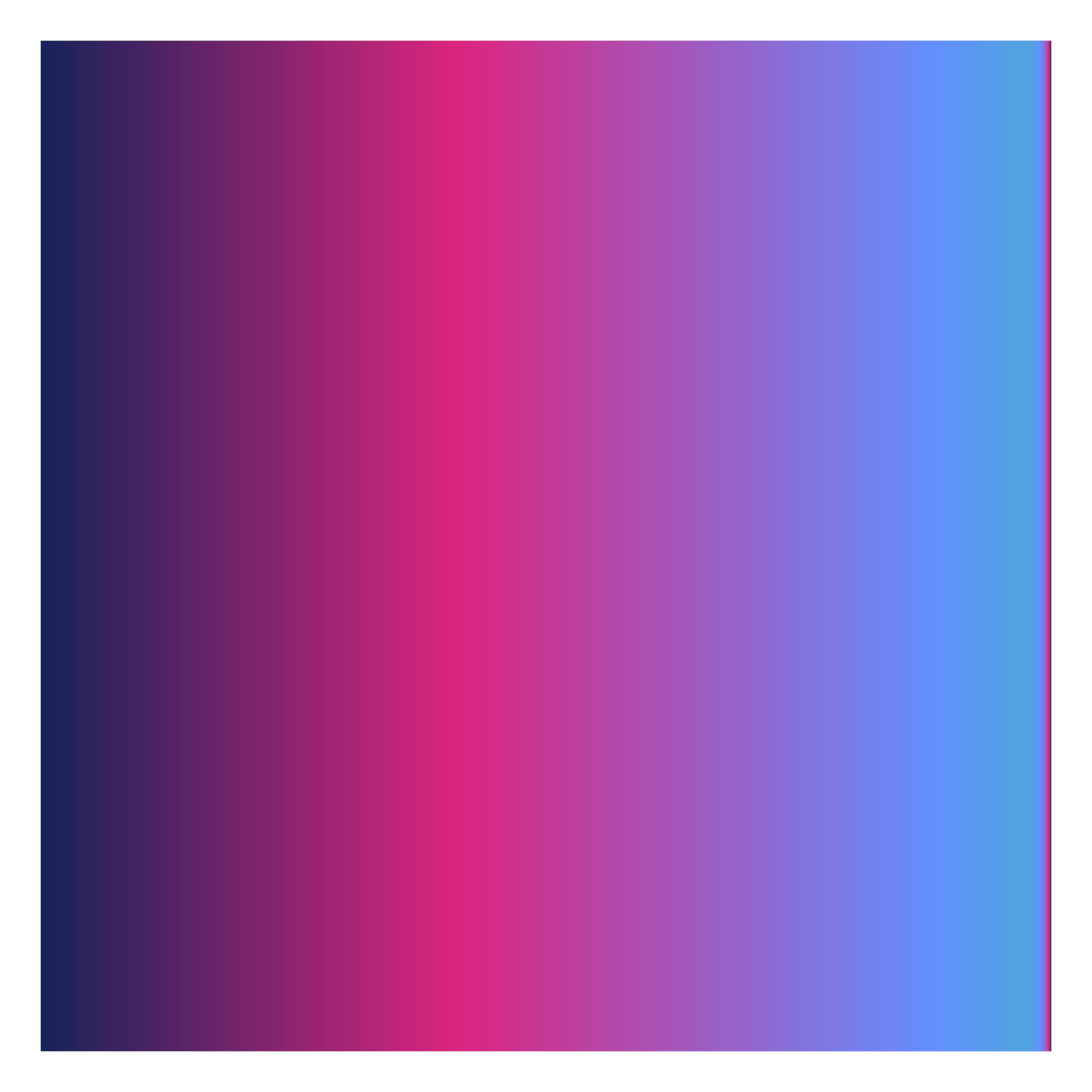}
        \includegraphics[scale=0.066]{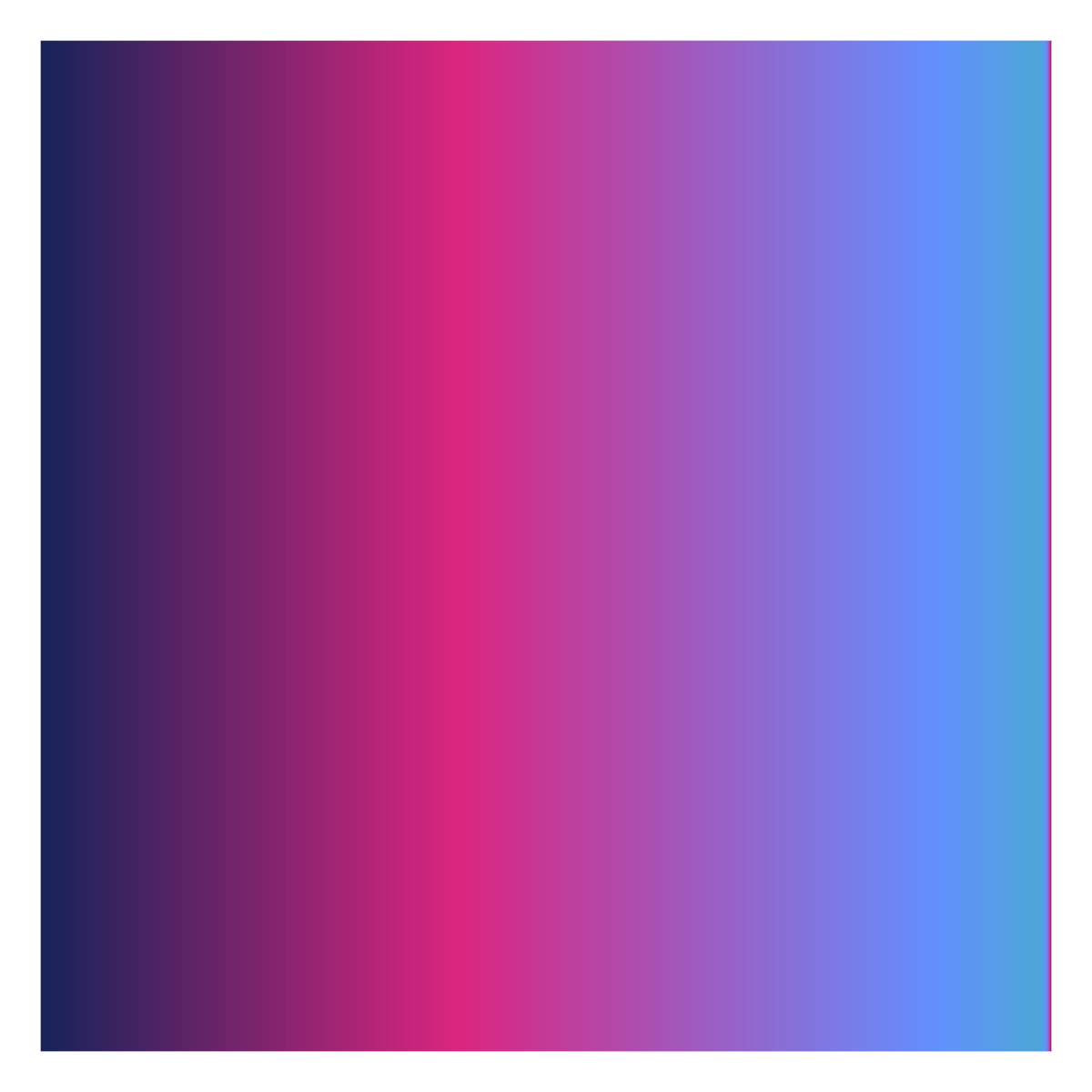}
        \includegraphics[scale=0.066]{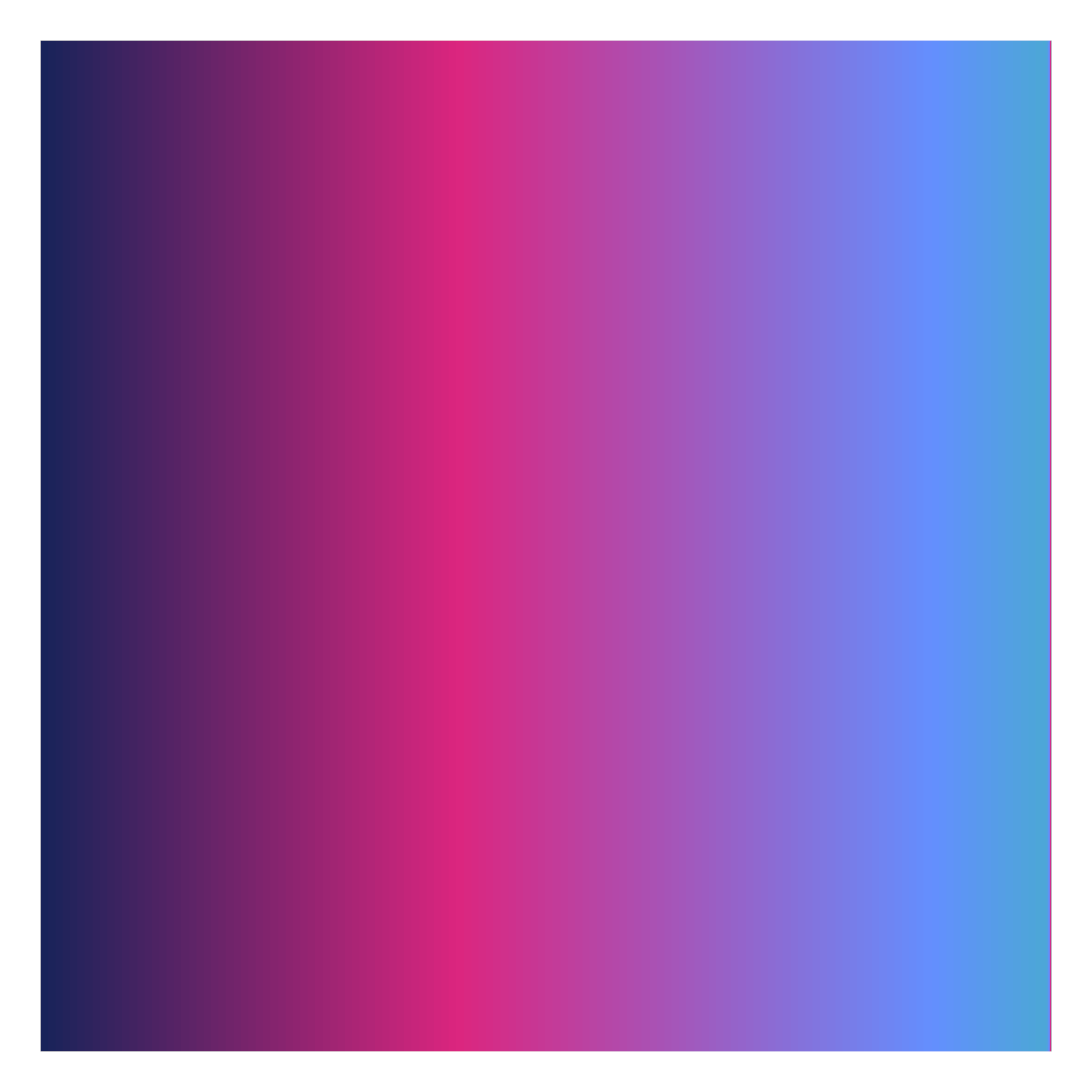}\hspace{3pt}\includegraphics[height=2.79cm]{figures/cw_sol_colorbar.png}

        \vspace{0.5cm}

        \includegraphics[scale=0.066]{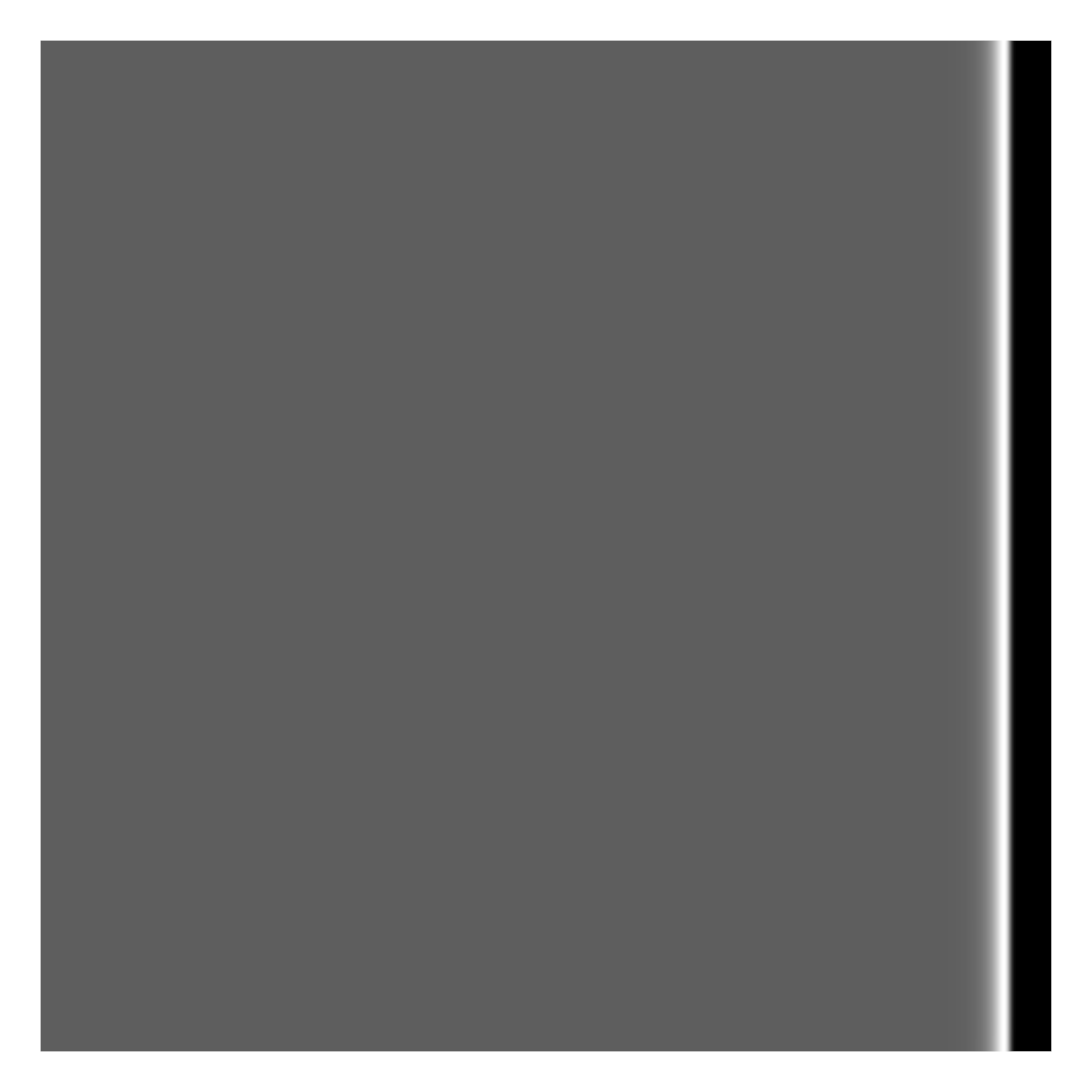}
        \includegraphics[scale=0.066]{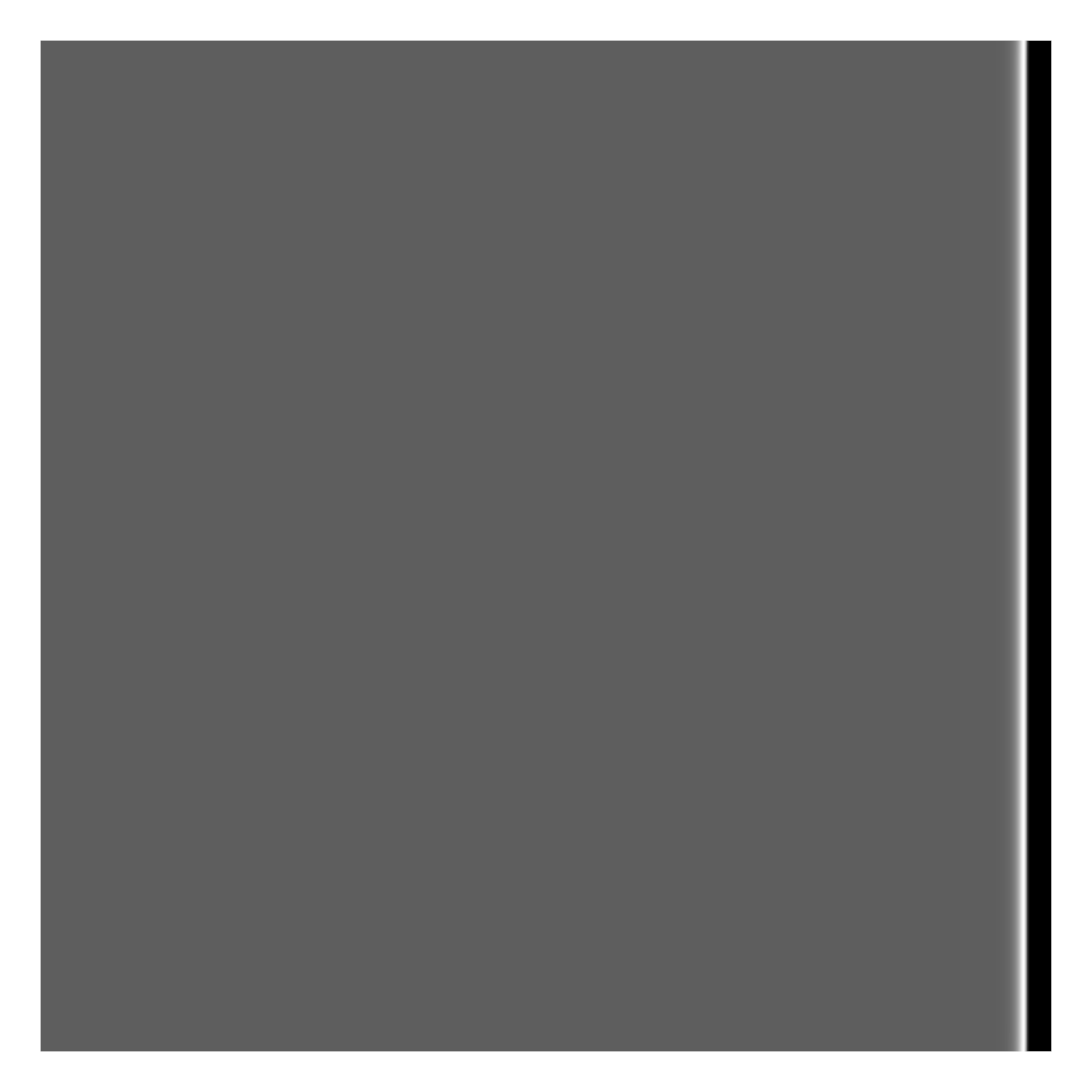}
        \includegraphics[scale=0.066]{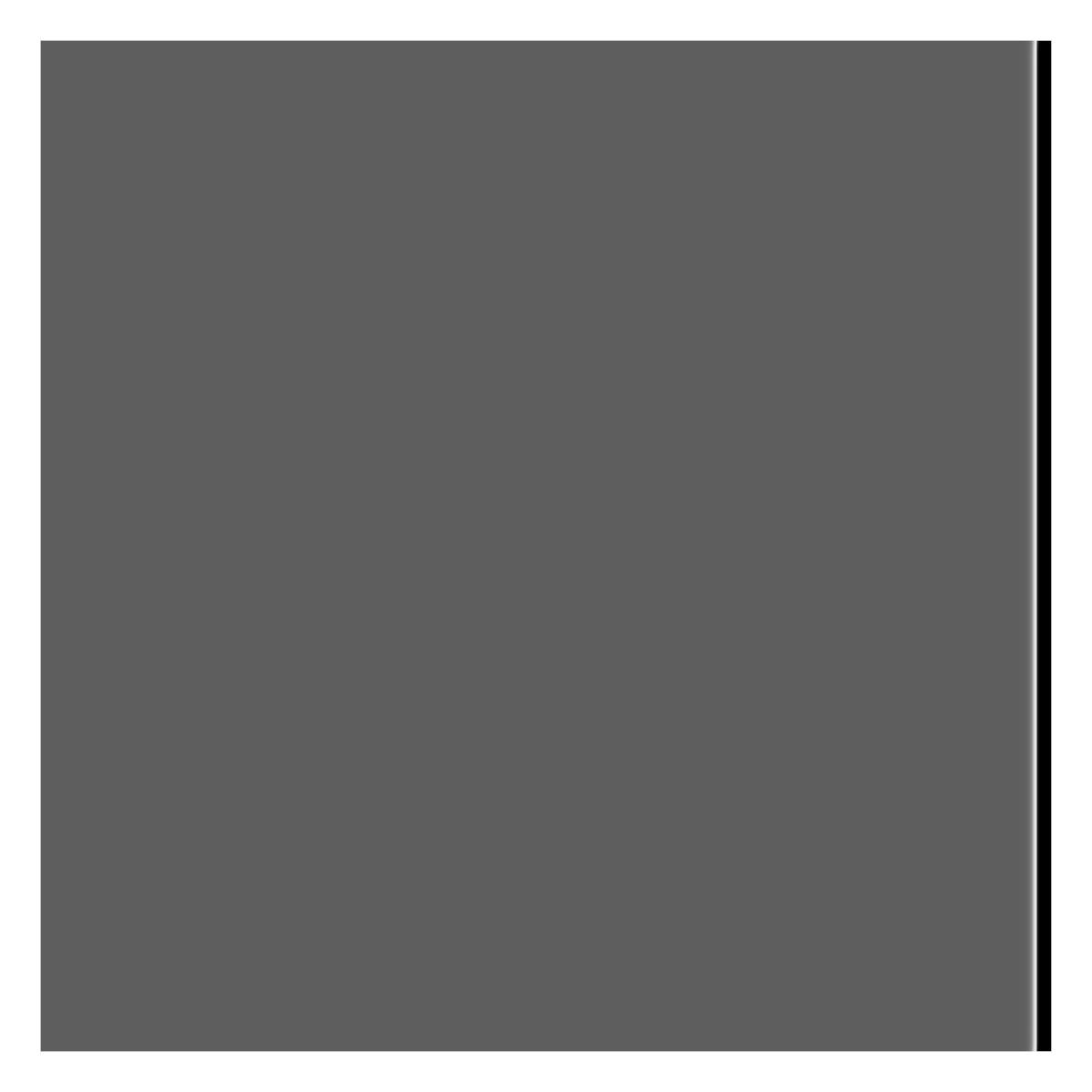}
        \includegraphics[scale=0.066]{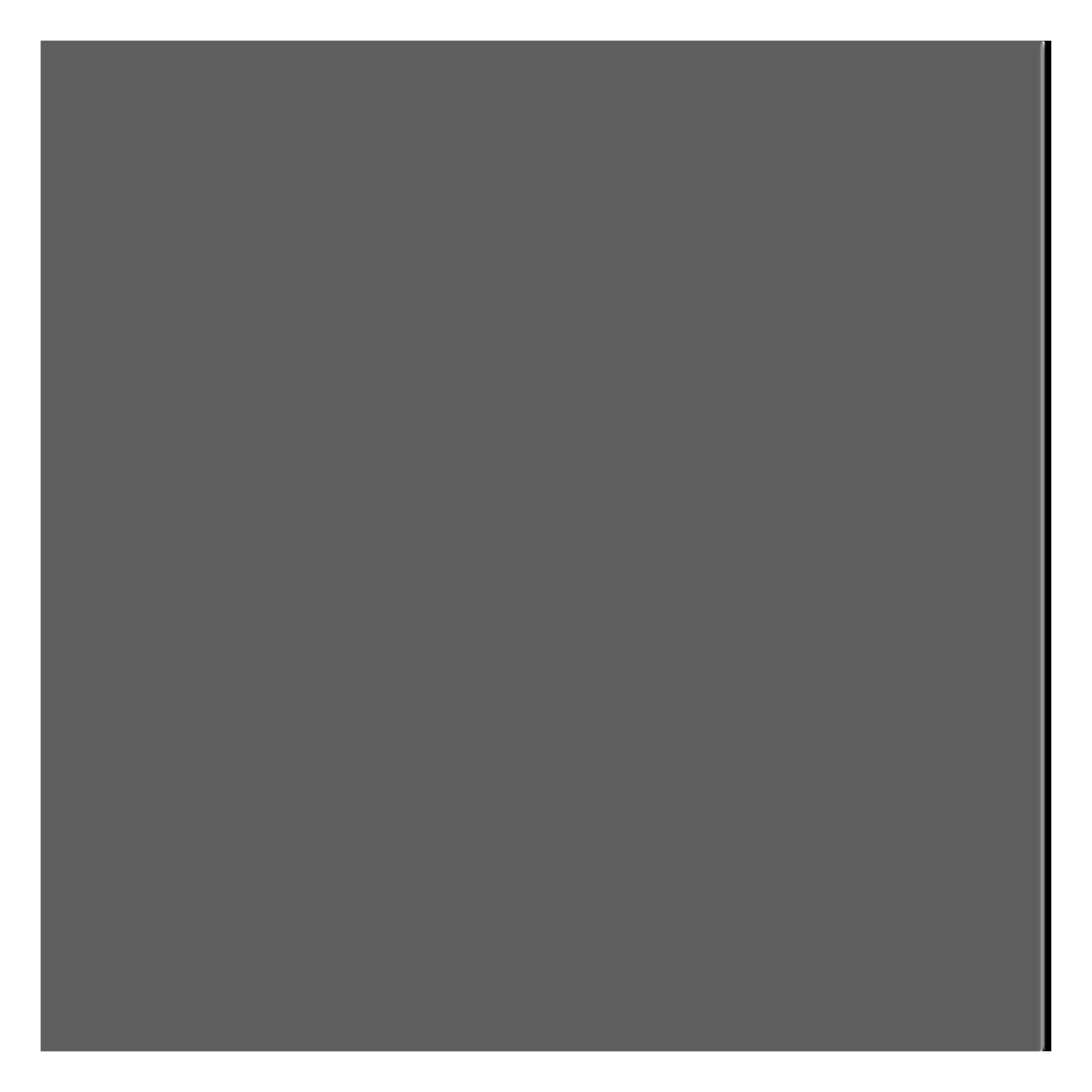}\hspace{3pt}\includegraphics[height=2.79cm]{figures/cw_div_colorbar.png}
        \caption{\bb{The discrete solution $u_h$ of the advection diffusion problem \eqref{eq:weakForm} with wind $\beta=(1,0)$ and mixed boundary conditions (homogeneous Dirichlet on the left and right boundaries, natural Neumann on the top and bottom), for different values of $\nu$ at the finest mesh size $512\times 512$, together with $\exp(-\lvert\nabla\cdot\beta u_h\rvert^2)$. Compared with the pure-Dirichlet problem of Figure~\ref{fig:reactionDiffusionAMG}, the top and bottom layers are absent, which is reflected in the lower iteration counts of Table~\ref{tab:reactionDiffusionAMG3} (see Remark~\ref{rmk:boundaryLayers}).}}
        \label{fig:reactionDiffusionAMG3}
    \end{figure}

    {To demonstrate that the methodology extends beyond constant winds, we conclude with the double-glazing benchmark of Elman, Silvester and Wathen \cite{Elman}. On $\Omega = (-1,1)^2$ the recirculating, divergence-free wind
    \begin{equation}
        \label{eq:glazingWind}
        \beta = \big(2y(1-x^2),\; -2x(1-y^2)\big)
    \end{equation}
    transports a scalar subject to a ``hot wall'' $u = 1$ on the right edge $x=1$ and $u=0$ on the remaining three edges, with no source term. As $\nu\to 0$ the solution develops both boundary and internal layers. The inhomogeneous Dirichlet data is imposed by lifting: writing $u = u_D + w$ with $u_D$ the hot-wall datum, we solve the homogeneous problem for the correction $w$, whose right-hand side is $-a(u_D, v_h)$.

    Because the wind is now solenoidal and genuinely two-dimensional, the recirculating flow has closed streamlines, along which the advection operator $\beta\cdot\nabla$ has a large near-kernel. A single GAMG V-cycle on the unprojected reaction diffusion operator \eqref{eq:discreteReactionDiffusion2}, which sufficed for the constant winds, degrades markedly as $\nu\to 0$: this streamline near-kernel is precisely where the sparse surrogate \eqref{eq:discreteReactionDiffusion2} differs from the projected normal PDE \eqref{eq:discreteReactionDiffusion}. The projected operator \eqref{eq:discreteReactionDiffusion} does not suffer from it: whenever $u_h$ is constant along streamlines, $\beta u_h$ is divergence free and hence $L^2$-orthogonal to gradients, so the projection $\Pi_\nabla$ annihilates it and the reaction term vanishes. Thus \eqref{eq:discreteReactionDiffusion} is the ideal preconditioner for the recirculating wind. The most prominent drawback of this approach is the computational cost, since the projection $\Pi_\nabla$ is global it requires the inverse Laplacian $S^{-1}$ which makes \eqref{eq:discreteReactionDiffusion} dense.

    Rather than approximate this dense operator by a sparse one, we apply it \emph{matrix-free}, giving it the same consistent SUPG streamline term $\delta N$ that stabilises $A$ (with $\delta = h/2|\beta|$, exactly as in the constant-wind experiments). Writing the preconditioner as
    \begin{equation}
        \label{eq:glazingMatFree}
        G = \nu S + C^{T} S^{-1} C + \delta N,
    \end{equation}
    with $S$ the stiffness matrix of $(\nabla\,\cdot\,,\nabla\,\cdot\,)$, $C$ the discretisation of $\nu^{-1/2}(\beta\,\cdot\,,\nabla\,\cdot\,)$, and $N$ the sparse streamline-diffusion matrix of $(\beta\cdot\nabla\,\cdot\,,\beta\cdot\nabla\,\cdot\,)$, each application of $G$ uses only the sparse matrices $S$, $C$, $C^{T}$, $N$ together with a single solve with $S$ for the action of the inverse Laplacian $S^{-1}$. This action is precisely the projection $\Pi_\nabla$ and can be evaluated by any fast Poisson solver, we here use a sparse Cholesky factorisation, though a multigrid solve would scale better. The streamline term $\delta N$ is a local finite element matrix, so it adds nothing to the cost of a matrix-vector product, yet it is essential: omitting it inflates the coarse-mesh iteration counts by an order of magnitude. We invert $G$ by an inner conjugate gradient iteration preconditioned by a single GAMG V-cycle on the sparse surrogate \eqref{eq:discreteReactionDiffusion2}, which is spectrally equivalent to \eqref{eq:discreteReactionDiffusion} (see \cref{apx:spectral_equivalence}). The solver is therefore nested: an outer CGNR iteration on the normal equation, preconditioned by $G^{-1}$, each application of which runs a short inner CG loop.

    Table~\ref{tab:glazing} reports the outer and inner iteration counts. The outer count is at most $\mathcal{O}(10)$ on the coarsest meshes and \emph{decreases} under refinement to just two iterations on the finest meshes for every viscosity. Although $G$ already carries the streamline term $\delta N$, this sparse term reproduces only the leading effect of the SUPG stabilisation on the normal operator $A^{T}TA$: writing $A = A_0 + \delta B$ with $A_0$ the unstabilised operator and $B$ the streamline term, one has $A^{T}TA = A_0^{T}TA_0 + \delta(A_0^{T}TB + BTA_0) + \delta^2 BTB$, whereas $G$ reproduces the ideal part $A_0^{T}TA_0$ plus the single term $\delta N$, the remaining $\mathcal{O}(\delta)$ are omitted. The residual mismatch is thus $\mathcal{O}(\delta) = \mathcal{O}(h)$, and since $\delta = h/2|\beta|$ vanishes under refinement, $G$ approaches the exact normal operator and $G^{-1}$ becomes an almost perfect preconditioner. This is why the outer count falls to two on the finest meshes, and why it is the coarse meshes, where $\delta$ is largest, that are the most demanding. The cost of the method is the total number of inner Poisson solves (shown in parentheses), which also decreases under refinement and levels off at a few tens of solves, so it is the coarse, under-resolved meshes that are the most expensive.}

    \begin{table}
        \centering
	\caption{\bb{Matrix-free projected preconditioner \eqref{eq:glazingMatFree} for the double-glazing benchmark on $(-1,1)^2$ with the recirculating wind \eqref{eq:glazingWind}. The exact projected operator $G=\nu S + C^{T}S^{-1}C + \delta N$, with the consistent SUPG streamline term ($\delta = h/2|\beta|$), is applied matrix-free---the inverse Laplacian $S^{-1}$ by a sparse Cholesky solve---and inverted by an inner CG preconditioned by a single GAMG V-cycle (twenty symmetric-SOR sweeps per level) on the sparse surrogate \eqref{eq:discreteReactionDiffusion2}; no coarse space or deflation is used. Each cell reports the number of outer CGNR iterations and, in parentheses, the total number of inner Poisson solves. The outer iteration is terminated when the unpreconditioned residual $\lVert A\vec{x}-\vec{b}\rVert_2$ of the original system fell below $10^{-5}$, the inner CG at relative tolerance $3\cdot 10^{-2}$.}}
        \label{tab:glazing}
        \begin{tabular}{c|c|c|c|c|c}
            \toprule
            $\nu$               & $32\times 32$ & $64\times 64$ & $128\times 128$ & $256\times 256$ & $512\times 512$ \\
            \midrule
            $1\cdot 10^{-2}$    & 8 (38)  & 5 (25)  & 3 (18)  & 3 (24) & 2 (19) \\
            $5\cdot 10^{-3}$    & 12 (99) & 6 (54)  & 4 (34)  & 3 (24) & 2 (23) \\
            $2.5\cdot 10^{-3}$  & 18 (221) & 8 (107) & 4 (56)  & 3 (29) & 2 (26) \\
            $1.25\cdot 10^{-3}$ & 27 (497) & 15 (381) & 5 (126) & 3 (49) & 2 (46) \\
            \bottomrule
        \end{tabular}
    \end{table}

\end{example}

We have shown that by choosing the Riesz map $\tau$ associated with the $H^1_0(\Omega)$ inner product, the normal PDE associated with the normal equation is a reaction diffusion equation which can be efficiently preconditioned via multigrid strategies.
We do not claim that the preconditioner obtained via normal preconditioning is better than a carefully designed preconditioner for the original linear system.
In fact, with ad-hoc strategies efficient geometric \cite{Elman} and algebraic \cite{air,air2,air3} preconditioners can be derived for the original linear system.
Yet, we believe that normal preconditioning offers a viable alternative to classical preconditioning strategies for the original linear system, since the self-adjoint and positive definite nature of the normal equation allows us to use a wide range of simpler preconditioning strategies, such as algebraic multigrid preconditioning, as shown in Example \ref{ex:H1InnerProduct}.
Thus, it should be regarded as a valid option for practitioners when dealing with PDEs that are difficult to precondition, for which fine tuning of classical preconditioners for the original linear system is not feasible.
\begin{remark}
    Although \eqref{eq:discreteReactionDiffusion} contains no advection term, the reaction contribution
    \begin{equation}
        \nu^{-1} (\Pi_{\nabla}(\beta u_h), \Pi_{\nabla} (\beta v_h)),
    \end{equation}
    is anisotropic, while the diffusion term in \eqref{eq:discreteReactionDiffusion} is isotropic. Thus, as $\nu \to 0$, the corresponding weak formulation may develop boundary layers along the entire boundary $\partial \Omega$.
    Consequently, a stable discretization of \eqref{eq:discreteReactionDiffusion} requires mesh refinement along the all $\partial \Omega$, in contrast to the advection diffusion operator, for which boundary layers occur only at the outflow boundary and refinement can thus be localized. This leads to a mismatch between the refinement requirements of the two operators when non-uniform meshes are adapted to the advection diffusion problem.
\end{remark}
\begin{remark}[Relations with LSQ-FEM]
    There is an intimate relation between \eqref{eq:normalPDE} and the Least SQuares Finite Element Method (LSQ-FEM) \cite{Bochev}.
    In fact, \eqref{eq:normalPDE} can be recast as a minimization problem, i.e
    \begin{equation}
        \vec{x} =  \underset{\vec{w}\in V_h}{\text{argmin}}\; J(\vec{w};\vec{b}), \qquad J(\vec{w};\vec{b})\coloneqq \norm{A\vec{w} - \vec{b}}_T^2.
    \end{equation}
    This effectively reduces the discretisation associated with the normal equation \eqref{eq:normalPDE} to a rather non-standard LSQ-FEM discretisation. We stress that this is non-standard because it is typical in LSQ-FEM to view the least square functional as a residual associated with the strong form of the PDE or with one of its first order reformulations.
    Yet, a key distinction between the strategy discussed here and a non-standard least square finite elements formulation lies in the fact that while a LSQ-FEM method construct a variational formulation by selecting appropriate norms, function spaces and residual, the strategy here seeks to identify the continuous PDE discretised by the operator \( A^T T A \) for different choices of inner product $T$ and pick the inner product $T$ resulting in the most convenient PDE to precondition.
\end{remark}

All codes used to 
generate the results presented in this section can be found in the repository \url{https://github.com/UZerbinati/normalPreconditioning}.
\section{Conclusions}
\label{sec:conclusion}
We have presented the idea of {normal preconditioner}, a notion of goodness of a preconditioner for the normal equation. As a consequence, we have illustrated the non-uniqueness of an ideal {normal preconditioner}. We believe that the large number of ideal {normal preconditioners} already makes a significant case as to why we should regard the problem of preconditioning the normal equation as feasible as preconditioning the original linear system.

Additionally, we gave practical examples of situations where the normal equation arising from PDEs can be preconditioned to obtain an iteration count comparable to classical preconditioners for the original linear system.  We named the general framework devised to construct such preconditioner: normal preconditioning.
The normal preconditioning framework can be summarised as follows:
\begin{enumerate}
    \item identify the correct class of Riesz maps $\tau:V^\prime\to V$ from the function space involved in the continuous PDE we have discretised;
    \item consider a PDE such that its discretisation $P$ is a good {normal preconditioner} for the normal equation;
    \item consider as preconditioner for the normal equation the matrix $P^TTP$.
\end{enumerate}
We investigated the advection diffusion equation discretized with both finite elements and finite differences.
In Example \ref{ex:H1InnerProduct}, following \textit{step 2}, we picked as $P$ a discretisation of the reaction diffusion operator, stabilised with the same streamline-diffusion term as the advection diffusion operator $A$, which we discretised with SUPG stabilisation.
We then showed that the preconditioner $P^TTP$ can easily be inverted via multigrid strategies, since it is a reaction-diffusion operator.
Lastly, we showed that this strategy yields a preconditioner that is mesh-independent for a sufficiently small value of $\nu$ and that requires a number of iterations comparable to the one obtained with classical preconditioners for the original linear system itself.

We would like to conclude by discussing some questions that we plan to investigate further in the future.
\begin{enumerate}
    \item \textbf{Mixed reformulations} of \eqref{eq:normalPDE} can be used to escape the fact that $A^TTA$ is a dense matrix. Specifically, we aim to study this problem from the point of view of its Schur complement \cite{Wathen}. Such reformulation has many similarities to the discontinuous Petrov--Galerkin and least squares finite element methods \cite{Demkowicz, Ita,Bochev}, which we plan to highlight in future work.
    \item \textbf{Different PDEs} can be tackled with the normal preconditioning framework, especially singularly perturbed problems. We aim to extend the presented ideas to the Oseen equation \cite{Elman}, the Helmholtz equation, and the \bb{Helmholtz}--Korteweg equation \cite{Farrell2}.
\end{enumerate}

\appendix
\section{Spectral equivalence}
\label{apx:spectral_equivalence}
In view of Theorem \ref{thm:conv} the convergence of the preconditioned CGNR method depends on the spectral condition number
\begin{equation}
    \kappa_T (AP^{-1}) = \frac{\sigma_{T,\textrm{max}}(AP^{-1})}{\sigma_{T,\textrm{min}}(AP^{-1})},
\end{equation}
where $\sigma_{\textrm{max}}(AP^{-1})$ and $\sigma_{T,\textrm{min}}(AP^{-1})$ are the maximal and minimal singular values of $AP^{-1}$ with respect to the inner product induced by $T$, as defined in \cref{sec:funcanal}.

While the theoretical spectral equivalence between \eqref{eq:discreteReactionDiffusion} and the normal  equation \eqref{eq:normalPDE} follow from the fact that both operators discretise the same continuous operator, the spectral equivalence between \eqref{eq:discreteReactionDiffusion2} and \eqref{eq:normalPDE} is not immediate.
We thus investigate the spectral equivalence of the proposed normal preconditioners numerically.
In particular, we compute approximations of the extremal eigenvalues associated with the generalized eigenvalue problem
\begin{equation}
    \label{eq:gen_eigenvalue_problem}
    A^{T} T A x = \sigma_T\, P^{T} T P x = \sigma_T\, G x,
\end{equation}
which characterizes the singular values of $AP^{-1}$ with respect to the inner product induced by $T$.
Notice in particular that we do not have access to the matrix $P$ associated with \eqref{eq:discreteReactionDiffusion2}, but we will construct the matrix $G$ associated with \eqref{eq:discreteReactionDiffusion2} via PETSc GAMG \cite{Adams}.

\begin{table}
    \centering
    \caption{Comparison of $\kappa_T$, the maximal and minimal eiganvlue of \eqref{eq:gen_eigenvalue_problem} for different values of $\nu$ and different mesh sizes, where $G$ is obtained by PETSc GAMG applied to \eqref{eq:discreteReactionDiffusion2} with homogeneous Dirichlet boundary conditions along the all $\partial \Omega$.
        The wind is fixed to $\beta = (1,0)$.}
    \label{tab:spectral_reactionDiffusionAMG}
    \begin{tabular}{c|c|c|c}
        \toprule
        $\nu$               & $32\times 32$ & $64\times 64$ & $128\times 128$ \\
        \midrule
        $1\cdot 10^{-2}$    & $\sigma_{T,\text{min}}:\;4.26\cdot 10^{-2}$ & $\sigma_{T,\text{min}}:\;3.78\cdot 10^{-2}$ & $\sigma_{T,\text{min}}:\;3.04\cdot 10^{-2}$ \\
                            & $\sigma_{T,\text{max}}:\;2.24\cdot 10^{0}$ & $\sigma_{T,\text{max}}:\;1.54\cdot 10^{0}$ & $\sigma_{T,\text{max}}:\;1.43\cdot 10^{0}$ \\
                            & $\kappa_T:\,5.27\cdot 10^{1}$ & $\kappa_T:\,4.07\cdot 10^{1}$ & $\kappa_T:\,4.70\cdot 10^{1}$ \\
        \midrule
        $5\cdot 10^{-3}$    & $\sigma_{T,\text{min}}:\;2.34\cdot 10^{-2}$ & $\sigma_{T,\text{min}}:\;2.07\cdot 10^{-2}$ & $\sigma_{T,\text{min}}:\;1.87\cdot 10^{-2}$ \\
                            & $\sigma_{T,\text{max}}:\;3.20\cdot 10^{0}$ & $\sigma_{T,\text{max}}:\;2.25\cdot 10^{0}$ & $\sigma_{T,\text{max}}:\;1.54\cdot 10^{0}$ \\
                            & $\kappa_T:\,1.37\cdot 10^{2}$ & $\kappa_T:\,1.09\cdot 10^{2}$ & $\kappa_T:\,8.25\cdot 10^{1}$ \\
        \midrule
        $2.5\cdot 10^{-3}$  & $\sigma_{T,\text{min}}:\;1.20\cdot 10^{-2}$ & $\sigma_{T,\text{min}}:\;1.12\cdot 10^{-2}$ & $\sigma_{T,\text{min}}:\;1.03\cdot 10^{-2}$ \\
                            & $\sigma_{T,\text{max}}:\;3.03\cdot 10^{0}$ & $\sigma_{T,\text{max}}:\;3.20\cdot 10^{0}$ & $\sigma_{T,\text{max}}:\;2.26\cdot 10^{0}$ \\
                            & $\kappa_T:\,2.51\cdot 10^{2}$ & $\kappa_T:\,2.85\cdot 10^{2}$ & $\kappa_T:\,2.20\cdot 10^{2}$ \\
        \midrule
        $1.25\cdot 10^{-3}$ & $\sigma_{T,\text{min}}:\;5.38\cdot 10^{-3}$ & $\sigma_{T,\text{min}}:\;5.81\cdot 10^{-3}$ & $\sigma_{T,\text{min}}:\;5.51\cdot 10^{-3}$ \\
                            & $\sigma_{T,\text{max}}:\;3.94\cdot 10^{0}$ & $\sigma_{T,\text{max}}:\;3.04\cdot 10^{0}$ & $\sigma_{T,\text{max}}:\;3.21\cdot 10^{0}$ \\
                            & $\kappa_T:\,7.32\cdot 10^{2}$ & $\kappa_T:\,5.24\cdot 10^{2}$ & $\kappa_T:\,5.82\cdot 10^{2}$ \\
        \bottomrule
    \end{tabular}
    \vspace{1cm}
\end{table}

\begin{table}
    \centering
    \caption{Comparison of $\kappa_T$, the maximal and minimal eiganvlue of \eqref{eq:gen_eigenvalue_problem} for different values of $\nu$ and different mesh sizes, where $G$ is obtained by PETSc GAMG applied to \eqref{eq:discreteReactionDiffusion2} with mixed boundary conditions.
        The wind is fixed to $\beta = (1,0)$.}
    \label{tab:spectral_reactionDiffusionAMG_mixed}
    \begin{tabular}{c|c|c|c}
        \toprule
        $\nu$               & $32\times 32$ & $64\times 64$ & $128\times 128$ \\
        \midrule
        $1\cdot 10^{-2}$    & $\sigma_{T,\text{min}}:\;4.35\cdot 10^{-2}$ & $\sigma_{T,\text{min}}:\;3.80\cdot 10^{-2}$ & $\sigma_{T,\text{min}}:\;3.06\cdot 10^{-2}$ \\
                            & $\sigma_{T,\text{max}}:\;4.13\cdot 10^{0}$ & $\sigma_{T,\text{max}}:\;2.69\cdot 10^{0}$ & $\sigma_{T,\text{max}}:\;2.09\cdot 10^{0}$ \\
                            & $\kappa_T:\,9.49\cdot 10^{1}$ & $\kappa_T:\,7.08\cdot 10^{1}$ & $\kappa_T:\,6.84\cdot 10^{1}$ \\
        \midrule
        $5\cdot 10^{-3}$    & $\sigma_{T,\text{min}}:\;2.35\cdot 10^{-2}$ & $\sigma_{T,\text{min}}:\;2.09\cdot 10^{-2}$ & $\sigma_{T,\text{min}}:\;1.87\cdot 10^{-2}$ \\
                            & $\sigma_{T,\text{max}}:\;5.68\cdot 10^{0}$ & $\sigma_{T,\text{max}}:\;4.14\cdot 10^{0}$ & $\sigma_{T,\text{max}}:\;2.69\cdot 10^{0}$ \\
                            & $\kappa_T:\,2.41\cdot 10^{2}$ & $\kappa_T:\,1.98\cdot 10^{2}$ & $\kappa_T:\,1.44\cdot 10^{2}$ \\
        \midrule
        $2.5\cdot 10^{-3}$  & $\sigma_{T,\text{min}}:\;1.20\cdot 10^{-2}$ & $\sigma_{T,\text{min}}:\;1.12\cdot 10^{-2}$ & $\sigma_{T,\text{min}}:\;1.03\cdot 10^{-2}$ \\
                            & $\sigma_{T,\text{max}}:\;6.31\cdot 10^{0}$ & $\sigma_{T,\text{max}}:\;5.72\cdot 10^{0}$ & $\sigma_{T,\text{max}}:\;4.14\cdot 10^{0}$ \\
                            & $\kappa_T:\,5.24\cdot 10^{2}$ & $\kappa_T:\,5.10\cdot 10^{2}$ & $\kappa_T:\,4.02\cdot 10^{2}$ \\
        \midrule
        $1.25\cdot 10^{-3}$ & $\sigma_{T,\text{min}}:\;5.37\cdot 10^{-3}$ & $\sigma_{T,\text{min}}:\;5.81\cdot 10^{-3}$ & $\sigma_{T,\text{min}}:\;5.51\cdot 10^{-3}$ \\
                            & $\sigma_{T,\text{max}}:\;6.15\cdot 10^{0}$ & $\sigma_{T,\text{max}}:\;6.35\cdot 10^{0}$ & $\sigma_{T,\text{max}}:\;5.71\cdot 10^{0}$ \\
                            & $\kappa_T:\,1.14\cdot 10^{3}$ & $\kappa_T:\,1.09\cdot 10^{3}$ & $\kappa_T:\,1.04\cdot 10^{3}$ \\
        \bottomrule
    \end{tabular}
    \vspace{1cm}
\end{table}

The resulting eigenvalue bounds provide a quantitative measure of the quality of the preconditioner and allow us to assess whether the spectral equivalence constants remain uniformly bounded with respect to mesh refinement and problem parameters.
\bb{The operators here are the SUPG-stabilised $A$ and the consistently stabilised surrogate \eqref{eq:discreteReactionDiffusion2}, so these bounds correspond exactly to the preconditioner used in Tables~\ref{tab:reactionDiffusionAMG}--\ref{tab:reactionDiffusionAMG3}. In particular, from Table \ref{tab:spectral_reactionDiffusionAMG} we observe that, for each fixed viscosity $\nu$, the condition number $\kappa_T(AP^{-1})$ is essentially independent of the mesh size: the extremal generalized eigenvalues $\sigma_{T,\text{min}}$ and $\sigma_{T,\text{max}}$ are both bounded away from $0$ and $\infty$ uniformly in the mesh, with $\sigma_{T,\text{min}}$ proportional to $\nu$, so that $\kappa_T = \mathcal{O}(\nu^{-1})$. Compared with a small fixed stabilisation, the consistent streamline term lifts $\sigma_{T,\text{max}}$ to an $\mathcal{O}(1)$ value somewhat above one and inflates $\kappa_T$ by a modest, mesh-independent factor; the eigenvalues nonetheless stay clustered, which---together with the stronger smoother---is why the iteration counts of Tables~\ref{tab:reactionDiffusionAMG}--\ref{tab:reactionDiffusionAMG3} remain low.
The condition number thus grows as the singular perturbation parameter $\nu$ decreases, but does not deteriorate under mesh refinement, showing that using the classical $L^2$ inner product as in \eqref{eq:discreteReactionDiffusion2}, inverted by GAMG, yields a normal preconditioner that is robust with respect to mesh refinement.
The same behaviour is observed for \eqref{eq:discreteReactionDiffusion2} with mixed homogeneous Dirichlet and Neumann boundary conditions in Table \ref{tab:spectral_reactionDiffusionAMG_mixed}.}

\section*{Acknowledgements}
The authors would like to thank P.~Brubeck, P.~E.~Farrell, J.~M\'{a}lek, C.~Parker, V.~Simoncini, A.~J.~Wathen, and S.~Zampini for the fruitful discussions and invaluable suggestions while writing this paper.
\bibliographystyle{siamplain}
\bibliography{references}
\end{document}

%% file: ex_shared.tex
\usepackage{lipsum}
\usepackage{amsfonts}
\usepackage{graphicx}
\usepackage{epstopdf}
\usepackage{mathtools}
\usepackage{algorithmic}
\usepackage[normalem]{ulem}
\ifpdf
  \DeclareGraphicsExtensions{.eps,.pdf,.png,.jpg}
\else
  \DeclareGraphicsExtensions{.eps}
\fi
\usepackage{tikz}
\usepackage{pgfplots,pgfplotstable}
\usepackage{longtable}
\usepgfplotslibrary{colorbrewer,groupplots}
\pgfplotsset{compat=1.18}

\pgfplotsset{
cycle list/Set1-5,
cycle multiindex* list={
mark list*\nextlist
Set1-5\nextlist
},
}

\pgfplotsset{
    discard if not/.style 2 args={
        x filter/.append code={
            \edef\tempa{\thisrow{#1}}
            \edef\tempb{#2}
            \ifx\tempa\tempb
            \else
                
            \fi
        }
    }
}

\newsiamremark{remark}{Remark}
\newsiamremark{hypothesis}{Hypothesis}
\crefname{hypothesis}{Hypothesis}{Hypotheses}
\newsiamthm{claim}{Claim}

\headers{
Preconditioned normal equations for PDEs 
}{L. Lazzarino, Y. Nakatsukasa, and U. Zerbinati}

\title{
Preconditioned normal equations 
for solving discretised partial differential equations
\thanks{
\funding{This work was funded by the Fog Research Institute under contract no.~FRI-454.
YN is supported by EPSRC grants EP/Y010086/1 and EP/Y030990/1.
LL is supported by the Ada Lovelace Centre
Programme at the Scientific Computing Department, STFC.
}}}

\author{Lorenzo Lazzarino \and Yuji Nakatsukasa \and Umberto Zerbinati\thanks{Mathematical Institute, University of Oxford, Oxford, UK  (\email{
lorenzo.lazzarino@maths.ox.ac.uk, 
nakatsukasa@maths.ox.ac.uk, 
zerbinati@maths.ox.ac.uk}).}}

\usepackage{amsopn}

\newcommand{\llnote}[1]{}
\newcommand{\bb}[1]{#1}

\newcommand{\ignore}[1]{}